\documentclass{amsart}
\usepackage{amsmath,amsthm,amsfonts,amssymb,amscd}

\newlength{\dinwidth}
\newlength{\dinmargin}

\newtheorem{theorem}{Theorem}[section]
\newtheorem{lemma}[theorem]{Lemma}
\newtheorem{corollary}[theorem]{Corollary}
\newtheorem{proposition}[theorem]{Proposition}

\theoremstyle{definition}

\newtheorem{example}[theorem]{Example}

\theoremstyle{remark}

\numberwithin{equation}{section}

\newcommand{\ep}{\varepsilon}

\begin{document}
\title[Weak mixing properties]{Weak mixing properties of vector
sequences}

\author{L\'aszl\'o Zsid\'o}
\address{Department of Mathematics, University of Rome ``Tor Vergata'',
Via della Ricerca Scientifica, 00133 Rome, Italy}
\email{zsido@axp.mat.uniroma2.it}
\date{May 15, 2005.}

\thanks{Supported by MIUR, GNAMPA-INDAM and EU}
\dedicatory{Dedicated to the memory of our colleague Gert K. Pedersen}

\begin{abstract}
Notions of weak and uniformly weak mixing (to zero) are defined for
bounded sequences in arbitrary Banach spaces. Uniformly weak mixing
for vector sequences is characterized by mean ergodic convergence
properties. This characterization turns out to be useful in the
study of multiple recurrence, where mixing properties of vector
sequences, which are not orbits of linear operators, are
investigated. For bounded sequences, which satisfy a certain
domination condition, it is shown that weak mixing to zero is
equivalent with uniformly weak mixing to zero.
\end{abstract}

\maketitle

\section{Introduction}
We recall that the \emph{upper density} $D^*(\mathcal A)$ and the
\emph{lower density} $D_*(\mathcal A)$ of some $\mathcal A \subset
\mathbb{N}:=\{\, 0\, ,\, 1\, ,\, 2\, ,\,\ldots\,\}$ are defined by
\[
D^*(\mathcal A):=\varlimsup_{n\to\infty}\, \frac 1{n+1}\,
\mathrm{card}\, \big(\mathcal A \cap [0,n]\big)\, ,\quad
D_*(\mathcal A):=\varliminf_{n\to\infty}\, \frac 1{n+1}\,
\mathrm{card}\, \big(\mathcal A \cap [0,n]\big)
\]
(see e.g. \cite{F2}, Chapter 3, \S 5 or \cite{Kr}, \S 2.3). If upper
and lower densities coincide then $\mathcal A$ is called having
\emph{density} $D(\mathcal A):=D^*(\mathcal A)=D_*(\mathcal A)\,$.

Clearly, for $\mathcal A\subset\mathbb{N}^*:=\mathbb{N}\,\backslash
\,\{ 0\} =\{\, 1\, ,\, 2\, ,\,\ldots\,\}$ we can use also the
formulas
\[
D^*(\mathcal A)=\varlimsup_{n\to\infty}\, \frac 1n\,
\mathrm{card}\, \big(\mathcal A \cap [1,n]\big)\, ,\quad
D_*(\mathcal A)=\varliminf_{n\to\infty}\, \frac 1n\,
\mathrm{card}\, \big(\mathcal A \cap [1,n]\big)\, .
\]

The upper (resp. lower) density of a sequence $(k_j)_{j\geq 1}$ in
$\mathbb{N}^*$ means the upper (resp. lower) density of the subset
$\{\, k_j\, ;\, j\geq 1\,\}$ of $\mathbb{N}^*\,$. It is easy to see
that the lower density of a strictly increasing $(k_j)_{j\geq 1}$
is $>0$ if and only if $\sup_{j\geq 1}\,\frac{k_j}j <+\infty\,$.

Let $X$ be a Banach space with dual space $X^*\,$. We shall say
that a sequence $(x_{\substack{{}\\ k}})_{\substack{{}\\ k\geq 1}}$
in $X$ is \emph{weakly mixing to zero} if
\begin{equation}\label{w.mix.}
\lim_{n\to\infty}\; \frac 1n\,
\sum_{k=1}^n\, |\langle x^*,x_{\substack{{}\\ k}}\rangle | = 0\;\,
\text{for all}\;\, x^*\in X^*\, ,
\end{equation}
and we shall say that it is \emph{uniformly weakly mixing to zero} if
\begin{equation}\label{unif.w.mix.}
\lim_{n\to\infty}\;
\sup\, \Big\{\, \frac 1n\, \sum\limits_{k=1}^n\, |\langle x^*,
x_{\substack{{}\\ k}}\rangle |\, ;\, x^*\in X^*\, ,\,\|x^*\|
\leq 1\,\Big\}= 0\, .
\end{equation}
A linear operator $U : X\to X$ is usually called \emph{weakly mixing
to zero at $x\in X$} if the orbit
$\big( U^k(x)\big)_{\substack{{}\\ k\geq 1}}$ is weakly mixing to
zero.

The following characterization of weak mixing to zero for power
bounded linear operators, which is a counterpart of the Blum-Hanson
theorem \cite{BlH} for weak mixing, was proved by L. K. Jones and
M. Lin \cite{JL} :
\begin{theorem}\label{orbit}
Let $U$ be a power bounded linear operator on a Banach space $X\,$,
$x\in X\,$, and $x_{\substack{{}\\ k}} = U^k(x)\, ,\, k\geq 1\,$.
Then the following conditions are equivalent $:$
\begin{itemize}
\item[(i)] The sequence
$(x_{\substack{{}\\ k}})_{\substack{{}\\ k\geq 1}}$ is
weakly mixing to zero.
\item[(j)] The sequence
$(x_{\substack{{}\\ k}})_{\substack{{}\\ k\geq 1}}$ is
uniformly weakly mixing to zero.
\item[(jj)] For every sequence $k_1<k_2<\, \ldots$ in $\mathbb{N}^*$
of lower density $>0\,$,
\[
\lim_{n\to\infty}\; \Big\|\,\frac 1n\, \sum_{k=1}^n\,
x_{\substack{{}\\ k_j}}\,\Big\| =0\; .\]
\end{itemize}
\end{theorem}

One main goal of this paper is to prove in the next section that
conditions (j) and (jj) in Theorem \ref{orbit} are equivalent for any
bounded sequence $(x_{\substack{{}\\ k}})_{\substack{{}\\ k\geq 1}}$
in the Banach space $X\,$, not only for the points of an orbit of
some power bounded linear operator on $X$ (Theorem \ref{unif.wm}).
Therefore, for any bounded sequence in a Banach space, uniformly weak
mixing to zero is equivalent with the mean ergodic convergence
property from (jj) (in particular, for bounded sequences in Hilbert
spaces, our notion of ``uniformly weak mixing to zero'' coincides
with the notion of ``weak mixing'' considered in \cite{BB}).

We notice that this result was used by C. Niculescu, A. Str\"oh
and L. Zsid\'o to prove that if $\Phi$ is a $*$-endomorphism of a
$C^*$-algebra $A\,$, leaving invariant a state $\varphi$ of $A\,$,
whose support in $A^{**}$ belongs to the centre of $A^{**}\,$, and
$\Phi$ is weakly mixing with respect to $\varphi\,$, then $\Phi$
is automatically weakly mixing of order $2$ (\cite{NSZ}, Theorem 1.3):
this is a partial extension to the non-commutative
$C^*$-dynamical systems of a classical result of H. Furstenberg,
according to which every weakly mixing measure preserving
transformation of a probability measure space is weakly mixing of
any order (\cite{F2}, Theorem 4.11).

For general bounded sequences in Banach spaces (or even in Hilbert
spaces), condition (i) in Theorem \ref{orbit} does not imply the
equivalent conditions (j) and (jj) (Examples \ref{ex1} and
\ref{ex2}). Nevertheless, we shall prove in Section \ref{c.sh-bdd}
that (i) implies (j) and (jj) provided that the sequence satisfies
some appropriate domination condition, called ``convex
shift-boundedness'', which of course holds if the sequence is an
orbit of some power bounded linear operator. Actually it will be
proved that if a convex shift-bounded sequence
$(x_{\substack{{}\\ k}})_{\substack{{}\\ k\geq 1}}$ in the Banach
space $X$ is weakly mixing to zero, then
\begin{equation}\label{B.unif.w.mix.}
\lim_{\substack{a,b\in\mathbb N^*\\ b-a\to\infty}}\;
\sup \Big\{\, \frac 1{b-a+1}\, \sum\limits_{k=a}^b\, |\langle x^*,
x_{\substack{{}\\ k}}\rangle |\, ;\, x^*\in X^*\, ,\,\|x^*\|
\leq 1\,\Big\}= 0
\end{equation}
(Theorem \ref{wm=unif.wm}). Its proof depends upon a structure
theorem for sets of natural numbers of non-zero upper Banach
density (Theorem \ref{structure}), which is of interest for itself.
We notice that if $(x_{\substack{{}\\ k}})_{\substack{{}\\ k\geq 1}}$
is an orbit of some power bounded linear operator on $X\,$, then
(\ref{B.unif.w.mix.}) is an immediate consequence of
(\ref{unif.w.mix.}).

Finally, in Section \ref{Cesaro} it will be shown that in uniformly
convex Banach spaces the above implication holds for sequences
which satisfy a condition weaker than convex shift-boundedness
(Theorem \ref{Cesaro-bdd}).

We notice that a short investigation of the ergodicity, that is
of the Cesaro norm-convergence to zero, of convex shift-bounded
sequences is postponed in an appendix.

\bigskip
\section{Uniformly weak mixing to zero}
A subset $\mathcal N$ of $\mathbb{N}^*$ is called \emph{relatively
dense} if there exists $L>0$ such that every interval of natural
numbers of lenght $\geq L$ contains some element of $\mathcal N\,$.
In this case holds clearly $D_*(\mathcal N)\geq \frac 1L\,$, so
relatively dense sets are of lower density $>0\,$.

A sequence $(k_j)_{j\geq 1}$ in $\mathbb{N}^*$ is called relatively
dense if the subset $\{\, k_j\, ;\, j\geq 1\,\}$ of $\mathbb{N}^*$
is relatively dense. It is easy to see that a strictly increasing
sequence $(k_j)_{j\geq 1}$ is relatively dense if and only if
$\sup_{j\geq 1}\, (k_{j+1} - k_j) <+\infty\,$.

The proof of the following lemma is immediate and we give it only
for the sake of completeness :

\begin{lemma}\label{sums}
For any sequence $(x_{\substack{{}\\ k}})_{\substack{{}\\ k\geq 1}}$
in a Banach space and any sequence $(k_j)_{j\geq 1}$ in $\mathbb{N}^*$
of lower density $>0$ we have
\[
\Big\|\,\frac 1n\, \sum_{j=1}^n\, x_{\substack{{}\\ k_j}}\,\Big\|
\longrightarrow 0 \;\iff\;
\Big\|\,\frac 1n\, \sum_{\substack{k\in\{k_1,k_2,\ldots\}\\ k\leq n}}
\, x_{\substack{{}\\ k}}\,\Big\|\longrightarrow 0\, .
\]
\end{lemma}

\begin{proof} For $\Longrightarrow\,$: with $n\geq k_1\,$, defining
$j(n)\in\mathbb N^*$ by $k_{j(n)}\leq n< k_{j(n)+1}\,$, we have
\begin{equation*}
\Big\|\,\frac 1n\, \sum_{\substack{k\in\{k_1,k_2,\ldots\}\\ k\leq n}}
\, x_{\substack{{}\\ k}}\,\Big\| =
\Big\|\,\frac 1n\, \sum_{j=1}^{j(n)}\, x_{\substack{{}\\ k_j}}\,\Big\|
= \underbrace{\frac{j(n)}n}_{\leq\, 1}\, \cdot\, \Big\|\,
\frac 1{j(n)}\, \sum_{j=1}^{j(n)}\, x_{\substack{{}\\ k_j}}\,\Big\|
\xrightarrow{\; n\to\infty\;} 0\, .
\end{equation*}
The converse implication $\Longleftarrow\,$ follows by using
\begin{equation*}
\Big\|\,\frac 1n\, \sum_{j=1}^n\, x_{\substack{{}\\ k_j}}\,\Big\| =
\Big\|\,\frac 1n\, \sum_{\substack{k\in\{k_1,k_2,\ldots\}\\ k\leq k_n}}
\, x_{\substack{{}\\ k}}\,\Big\| = \frac{k_n}n \cdot\, \Big\|\,
\frac 1{k_n}\, \sum_{\substack{k\in\{k_1,k_2,\ldots\}\\ k\leq k_n}}
\, x_{\substack{{}\\ k}}\,\Big\|\, .
\end{equation*}

\end{proof}

The next lemma is the main ingredient in the proof of the main result
of the section :
\begin{lemma}\label{key}
Let $\Omega$ be a compact Hausdorff topological space, and $f_1\, ,\,
f_2\, ,\, \ldots\,$ \\ continuous complex functions on $\Omega$ of
uniform norm $\| f_{\substack{{}\\ k}}\|_{\substack{{}\\ \infty}}
\leq 1\,$. If
\[
\Big\|\,\frac 1n\, \sum_{j=1}^n\, f_{\substack{{}\\ k_j}}\,
\Big\|_\infty \longrightarrow 0\;\text{ for every relatively dense }\;
(k_j)_{j\geq 1}\subset\mathbb{N}^*
\]
then
\[
\Big\|\,\frac 1n\, \sum_{k=1}^n\, |f_{\substack{{}\\ k}}|\,
\Big\|_\infty \longrightarrow 0 \, .
\]
\end{lemma}

\begin{proof} Without loss of generality we can assume that the
functions $f_{\substack{{}\\ k}}$ are real. Furthermore, since
$|f_{\substack{{}\\ k}}| = 2\,{f_{\substack{{}\\ k}}}^+ -
f_{\substack{{}\\ k}}\,$, it is enough to prove that
\[
\Big\|\,\frac 1n\, \sum_{k=1}^n\, {f_{\substack{{}\\ k}}}^+\,
\Big\|_\infty \longrightarrow 0 \, .
\]
Let us assume the contrary, that is the existence of some $\ep_o >0$
for which
\[
\mathcal J := \Big\{\, n\geq 1\, ;\, \Big\|\,\frac 1n\, \sum_{k=1}^n\,
{f_{\substack{{}\\ k}}}^+\,\Big\|_\infty \geq\ep_o\,\Big\}\;
\text{ is infinite.}
\]
For every $n\in\mathcal J$ there exists $\omega_n \in\Omega$ such that
\[
\text{the cardinality of }\; \mathcal N_n := \{\, 1\leq k\leq n\, ;\,
{f_{\substack{{}\\ k}}}^+(\omega_n)\geq \frac{\ep_o}2 \,\}\;
\text{ is }\; \geq\frac{n\,\ep_o}2\, .
\]
Indeed, if $\omega_n \in\Omega$ is chosen such that
\[
\frac 1n\, \sum_{k=1}^n\, {f_{\substack{{}\\ k}}}^+(\omega_n ) =
\Big\|\,\frac 1n\, \sum_{k=1}^n\, {f_{\substack{{}\\ k}}}^+\,
\Big\|_\infty \geq\ep_o
\]
then
\begin{align*}
\ep_o &\leq \frac 1n\,  \Big( \sum_{k\in\mathcal N_n}\,
{f_{\substack{{}\\ k}}}^+(\omega_n ) + \sum_{\substack{1\leq k\leq n \\
k\notin\mathcal N_n}}\, {f_{\substack{{}\\ k}}}^+(\omega_n )\Big)\leq \\
&\leq \frac 1n\, \Big(\mathrm{card}\, (\mathcal N_n ) +
\frac{\ep_o}2\, \big(n-\mathrm{card}\, (\mathcal N_n )\big)\Big)\leq \\
&\leq \frac 1n\,\mathrm{card}\, (\mathcal N_n ) +\frac{\ep_o}2\, .
\end{align*}

Denoting now the least element of $\mathcal J$ by $k_1\,$, we can
construct recursively a sequence $k_1<k_2<\, \ldots\,$ in $\mathcal J$
such that
\[
\text{the cardinality of }\; {\mathcal N_{\substack{{}\\ k_{j+1}}}}' :=
\{\, k\in\mathcal N_{\substack{{}\\ k_{j+1}}}\, ;\, k > k_j \,\}\;
\text{ is }\; \geq\frac{k_{j+1}\cdot\ep_o}4\, ,\quad j\geq 1\, .
\]
For it is enough to choose $k_{j+1}\geq \frac{4\, k_j}{\ep_o}\,$,
because then
\[
\mathrm{card}\, ({\mathcal N_{\substack{{}\\ k_{j+1}}}}') \geq
\mathrm{card}\, (\mathcal N_{\substack{{}\\ k_{j+1}}}) - k_j \geq
\frac{k_{j+1}\cdot\ep_o}2 - k_j \geq\frac{k_{j+1}\cdot\ep_o}4\, .
\]
Putting
\[
\mathcal N := \bigcup_{j\geq 2}\, {\mathcal N_{\substack{{}\\ k_j}}}'\, ,
\]
we have for every $j\geq 2$
\[
\mathcal N \cap (k_{j-1},k_j] = {\mathcal N_{\substack{{}\\ k_j}}}'
\subset \mathcal N_{\substack{{}\\ k_j}}\, ,
\]
in particular,
\begin{align*}
k\in\mathcal N\, ,\, k_{j-1}<k\leq k_j &\;\Longrightarrow\;
{f_{\substack{{}\\ k}}}^+ (\omega_{\substack{{}\\ k_j}})\geq
\frac{\ep_o}2\\
&\;\Longrightarrow\; f_{\substack{{}\\ k}} (\omega_{\substack{{}\\ k_j}})
= {f_{\substack{{}\\ k}}}^+ (\omega_{\substack{{}\\ k_j}})\geq
\frac{\ep_o}2\, .
\end{align*}

Let us choose some integer $p\geq \frac{16}{{\ep_o}^2}\,$. Since
\[
{\mathcal N}^{(p)} := \mathcal N \cup \{\, p\, ,\, 2\, p\, ,\, 3\, p\,
\ldots\,\}\subset\mathbb{N}^*
\]
is relatively dense, by the assumption on the functions
$f_{\substack{{}\\ k}}$ and by Lemma \ref{sums} there exists
$m_o\geq 1$ such that
\[
m\geq m_o\;\Longrightarrow\;\Big\|\,\frac 1m\,
\sum_{\substack{k\in{\mathcal N}^{(p)}\\ k\leq m}}\,
f_{\substack{{}\\ k}}\,\Big\|_\infty\leq \frac{{\ep_o}^2}{34}\, .
\]
Then we get for any $j\geq 2$ with $k_{j-1}\geq m_o$
\begin{align*}
\frac{{\ep_o}^2}{17} = 2\,\frac{{\ep_o}^2}{34} &\geq\Big\|\,\frac 1{k_j}\,
\sum_{\substack{k\in{\mathcal N}^{(p)}\\ k\leq k_j}}\,
f_{\substack{{}\\ k}}\,\Big\|_\infty + \Big\|\,\frac 1{k_{j-1}}\,
\sum_{\substack{k\in{\mathcal N}^{(p)}\\ k\leq k_{j-1}}}\,
f_{\substack{{}\\ k}}\,\Big\|_\infty \geq \\
&\geq \Big\|\,\frac 1{k_j}\,
\sum_{\substack{k\in{\mathcal N}^{(p)}\\ k_{j-1}<k\leq k_j}}\,
f_{\substack{{}\\ k}}\,\Big\|_\infty \geq \Big|\,\frac 1{k_j}\,
\sum_{\substack{k\in{\mathcal N}^{(p)}\\ k_{j-1}<k\leq k_j}}\,
f_{\substack{{}\\ k}}(\omega_{\substack{{}\\ k_j}})\,\Big| \geq \\
&\geq \Big|\,\frac 1{k_j}\,
\sum_{\substack{k\in{\mathcal N}\\ k_{j-1}<k\leq k_j}}\,
f_{\substack{{}\\ k}}(\omega_{\substack{{}\\ k_j}})\,\Big| -
\Big|\,\frac 1{k_j}\,
\sum_{\substack{k_{j-1}<k\leq k_j\\ k\;\text{multiple of}\; p}}\,
f_{\substack{{}\\ k}}(\omega_{\substack{{}\\ k_j}})\,\Big| \geq \\
&\geq \frac 1{k_j}\cdot\frac{\ep_o}2\cdot \mathrm{card}\,
({\mathcal N_{\substack{{}\\ k_j}}}') -
\frac 1{k_j}\cdot \mathrm{card}\, (\{\, 1\leq k\leq k_j\, ;\,
k\;\text{multiple of}\; p\,\}) \geq \\
&\geq \frac 1{k_j}\cdot\frac{\ep_o}2\cdot \frac{k_j\cdot\ep_o}4
- \frac 1{k_j}\cdot\frac{k_j}p = \frac{{\ep_o}^2}8 - \frac 1p \geq
\frac{{\ep_o}^2}{16}\, ,
\end{align*}
which is absurde.

\end{proof}

Now we can characterize uniformly weak mixing to zero for 
bounded sequences in Banach spaces by mean ergodic convergence
properties :

\begin{theorem}[Mean ergodic description of uniformly weak
mixing]\label{unif.wm}
For a bounded sequence
$(x_{\substack{{}\\ k}})_{\substack{{}\\ k\geq 1}}$ in a
Banach space $X\,$, the following conditions are equivalent $:$
\begin{itemize}
\item[(j)] $(x_{\substack{{}\\ k}})_{\substack{{}\\ k\geq 1}}$ is
uniformly weakly mixing to zero, that is

\begin{center} $\displaystyle \lim_{n\to\infty}\;
\sup\, \Big\{\, \frac 1n\, \sum\limits_{k=1}^n\, |\langle x^*,
x_{\substack{{}\\ k}}\rangle |\, ;\, x^*\in X^*\, ,\,\|x^*\|
\leq 1\,\Big\}= 0\,$.
\end{center}
\item[(jj)] For every sequence $k_1<k_2<\, \ldots$ in $\mathbb{N}^*$
of lower density $>0\,$,

\begin{center} $\displaystyle \lim_{n\to\infty}\; \Big\|\,\frac 1n\,
\sum_{j=1}^n\, x_{\substack{{}\\ k_j}}\,\Big\| =0\;$.
\end{center}
\item[(jjj)] For every relatively dense sequence $k_1<k_2<\, \ldots$
in $\mathbb{N}^*\,$,

\begin{center} $\displaystyle \lim_{n\to\infty}\; \Big\|\,\frac 1n\,
\sum_{j=1}^n\, x_{\substack{{}\\ k_j}}\,\Big\| =0\;$.
\end{center}
\end{itemize}
\end{theorem}

\begin{proof} Implication (j)$\,\Rightarrow\,$(jj) follows
immediately from Lemma \ref{sums} and (jj)$\,\Rightarrow\,$(jjj)
is trivial.

For (jjj)$\,\Rightarrow\,$(j) we recall that the closed unit ball
$B_{\substack{{}\\ X^*}}$ of $X^*$ is weak${}^*$-compact and the
evaluation functions $f_x : B_{\substack{{}\\ X^*}}\ni x^*\mapsto\,
\langle x^*,x\rangle\, ,\, x\in X$ are weak${}^*$-continuous. Since
\begin{itemize}
\item[(j)] means $\displaystyle\, \frac 1n\, \sum_{k=1}^n\,
|f_{x_{\substack{{}\\ k}}}| \xrightarrow{\;\text{uniformly}\;} 0$ and
\item[(jjj)] means that, for every relatively dense sequence
$k_1<k_2<\, \ldots\,$ in $\mathbb{N}^*\,$,
\[
\frac 1n\, \sum_{k=1}^n\, f_{x_{\substack{{}\\ k_j}}}
\xrightarrow{\;\text{uniformly}\;} 0\, , \]
\end{itemize}
implication (jjj)$\,\Rightarrow\,$(j) follows from Lemma \ref{key}.

\end{proof}

Theorem \ref{unif.wm} yields a similar characterization of weak
mixing to zero :

\begin{corollary}\label{wm}
For a bounded sequence
$(x_{\substack{{}\\ k}})_{\substack{{}\\ k\geq 1}}$ in a
Banach space $X$ and $x^*\in X^*$, the following conditions are
equivalent $:$
\smallskip

\noindent\hspace{0.675 cm}{\rm (i)}${}_{\substack{{}\\ x^*}}$
\hspace{3.37 cm}$\displaystyle \lim_{n\to\infty}\; \frac 1n\,
\sum_{k=1}^n\, |\langle x^*,x_{\substack{{}\\ k}}\rangle | = 0\;$.
\begin{itemize}
\item[(ii)]$\!\!\!{}_{\substack{{}\\ x^*}}$ For every sequence
$k_1<k_2<\, \ldots$ in $\mathbb{N}^*$ of lower density $>0\,$,

\begin{center} $\displaystyle \lim_{n\to\infty}\;\Big\langle x^*,
\frac 1n\,\sum_{j=1}^n\, x_{\substack{{}\\ k_j}}\Big\rangle =0\;$.
\end{center}
\item[(iii)]$\!\!\!{}_{\substack{{}\\ x^*}}$ For every relatively
dense sequence $k_1<k_2<\, \ldots$ in $\mathbb{N}^*\,$,

\begin{center} $\displaystyle \lim_{n\to\infty}\;\Big\langle x^*,
\frac 1n\,\sum_{j=1}^n\, x_{\substack{{}\\ k_j}}\Big\rangle =0\;$.
\end{center}
\end{itemize}
\end{corollary}

\begin{proof}
We have just to apply Theorem \ref{unif.wm} to the bounded scalar
sequence
\smallskip

\noindent\hspace{5.24 cm}$\big(\langle x^*,x_{\substack{{}\\ k}}
\rangle\big)_{\substack{{}\\ k\geq 1}}\,$.

\end{proof}

\bigskip
\section{Comparison of weak and uniformly weak mixing to zero}

Let us first give an example of a bounded sequence in the Banach
space $C\big( [0,1]\big)$ of all continuous functions on
$[0,1]\,$, which satisfies (i) but not (j) in Theorem \ref{orbit}.
$\|\cdot\|_{\substack{{}\\ \infty}}$ will stand for the uniform
norm on $C\big( [0,1]\big)$ and supp$\, (f)$ will denote the
support of $f\in C\big( [0,1]\big)\,$.

\begin{example}\label{ex1}
{\it Let $1=n_{\substack{{}\\ 1}}<n_{\substack{{}\\ 2}}<\, \ldots\,$
be a sequence in $\mathbb{N}^*$ such that
\begin{center}
$\displaystyle \frac{n_j-1}{n_{j+1}-1}\leq\frac 12\, ,\qquad j\geq 1$
\end{center}
$($for example, $n_{\substack{{}\\ 1}}=1\,$, $n_{\substack{{}\\ 2}}=2$
and $n_{j+1}=2\, n_j -1\,$ for $\, j\geq 2\, )\,$,
\begin{center}
$1>t_{\substack{{}\\ 1}}>t_{\substack{{}\\ 2}}>\,\ldots\, >0\, ,
\quad t_j\longrightarrow 0$
\end{center}
real numbers, and $g_j : [0,1]\longrightarrow [0,1]\,$, $j\geq 1\,$,
continuous functions such that
\medskip

\noindent\hspace{2.34 cm}$\text{\rm supp}\, (g_j)\subset [t_{j+1},t_j]
\,$ and $\,\| g_j\|_{\substack{{}\\ \infty}} =1\,$ for all
$\, j\geq 1\, .$
\smallskip

\noindent If we set
\smallskip

\noindent\hspace{4.16 cm}$f_{\substack{{}\\ k}} =
g_{\substack{{}\\ j}}\,$ for $\, n_j\leq k<n_{j+1}\, ,$
\medskip

\noindent then $(f_{\substack{{}\\ k}})_{\substack{{}\\ k\geq 1}}$
is a bounded sequence in $\, C\big( [0,1]\big)$, which is weakly
convergent to zero, and so is weakly mixing to zero, but which is
not uniformly weakly mixing to zero.}
\end{example}

\begin{proof} Since $0\leq f_{\substack{{}\\ k}}\leq 1\,$ for every
$\, k\geq 1\,$, according to the Riesz representation theorem and
the Lebesgue dominated convergence theorem, the weak convergence
of $(f_{\substack{{}\\ k}})_{\substack{{}\\ k\geq 1}}$ to zero is
equivalent to
\begin{equation}\label{p.w.}
f_{\substack{{}\\ k}}\xrightarrow{\;\text{pointwise}\;} 0\, ,
\end{equation}
while (\ref{unif.w.mix.}) for
$(f_{\substack{{}\\ k}})_{\substack{{}\\ k\geq 1}}$
is equivalent with
\begin{equation}\label{unif.}
\frac 1n\, \sum_{k=1}^n\, f_{\substack{{}\\ k}}
\xrightarrow{\;\text{uniformly}\;} 0\, .
\end{equation}

For (\ref{p.w.}) let $t\in [0,1]$ be arbitrary. If $t=0$ then
(\ref{p.w.}) holds obviously because $f_{\substack{{}\\ k}} (0)
=0\,$ for all $ k\geq 1\,$. On the other hand, if $0<t\leq 1$
then there exists some $j\geq 1$ with $t_j<t$ and so
\smallskip

\noindent\hspace{4.575 cm}$f_{\substack{{}\\ k}} (t) =0\, ,\qquad
n\geq n_j\, .$
\medskip

%
%

Now, by the positivity of the functions $g_{\substack{{}\\ j}}$
and $f_{\substack{{}\\ k}}\,$, we have for every $j\geq 1\,$:
\begin{align*}
\frac 1{n_{j+1}-1}\sum_{k=1}^{n_{j+1}-1}f_{\substack{{}\\ k}}
&\geq\frac 1{n_{j+1}-1}\sum_{k=n_j}^{n_{j+1}-1}f_{\substack{{}\\ k}}
=\frac{n_{j+1}-n_j}{n_{j+1}-1}\, g_{\substack{{}\\ j}} =
\Big( 1-\frac{n_j-1}{n_{j+1}-1}\Big) g_{\substack{{}\\ j}} \\
&\geq\frac 12\, g_{\substack{{}\\ j}}\, .
\end{align*}
Consequently
\smallskip

\noindent\hspace{2.63 cm}$\displaystyle \Big\|\,\frac 1{n_{j+1}-1}
\sum_{k=1}^{n_{j+1}-1}f_{\substack{{}\\ k}}\Big\|_\infty \geq
\frac 12\, \| g_{\substack{{}\\ j}}\|_{\substack{{}\\ \infty}}
=\frac 12\,$ for all $\, j\geq 1$
\medskip

\noindent and so (\ref{unif.}) does not hold.

\end{proof}

A similar counterexample can be given also in the Hilbert space
$L^2\big( [0,1]\big)\,$, whose inner product and norm will be denoted
by $(\,\cdot\, |\,\cdot\, )$ and $\|\cdot\|_{\substack{{}\\2}}\,$,
respectively :

\begin{example}\label{ex2}
{\it Let $1=n_{\substack{{}\\ 1}}<n_{\substack{{}\\ 2}}<\, \ldots\,$
be a sequence in $\mathbb{N}^*$ such that
\begin{center}
$\displaystyle \frac{n_j-1}{n_{j+1}-1}\leq\frac 12\, ,\qquad j\geq 1$
\end{center}
$($for example, $n_{\substack{{}\\ 1}}=1\,$, $n_{\substack{{}\\ 2}}=2$
and $n_{j+1}=2\, n_j -1\,$ for $\, j\geq 2\, )\,$,
\begin{center}
$1>t_{\substack{{}\\ 1}}>t_{\substack{{}\\ 2}}>\,\ldots\, >0\, ,
\quad t_j\longrightarrow 0$
\end{center}
real numbers, and $g_j : [0,1]\longrightarrow [0,+\infty )\,$,
$j\geq 1\,$, continuous functions such that
\medskip

\noindent\hspace{2.34 cm}$\text{\rm supp}\, (g_j)\subset [t_{j+1},t_j]
\,$ and $\,\| g_j\|_{\substack{{}\\2}} =1\,$ for all $\, j\geq 1\, .$
\smallskip

\noindent If we set
\smallskip

\noindent\hspace{4.16 cm}$f_{\substack{{}\\ k}} =
g_{\substack{{}\\ j}}\,$ for $\, n_j\leq k<n_{j+1}\, ,$
\medskip

\noindent then $(f_{\substack{{}\\ k}})_{\substack{{}\\ k\geq 1}}$
is a bounded sequence in $L^2\big( [0,1]\big)$, which is weakly
convergent to zero, and so is weakly mixing to zero, but which is
not uniformly weakly mixing to zero.}
\end{example}

\begin{proof} Since the functions $g_j$ are mutually orthogonal,
by the Bessel inequality we have for every $f\in L^2\big( [0,1]
\big)\,$:
\smallskip

\noindent\hspace{2.527 cm}$\displaystyle \sum_{j=1}^\infty\,
|(g_j | f)|^2\leq\| f\|_{\substack{{}\\ 2}}^{\, 2}<+\infty\, ,\,$
hence $\; (g_j | f)\longrightarrow 0\,$.
\smallskip

\noindent Therefore $f_{\substack{{}\\ k}}
\xrightarrow{\;\text{weakly}\;} 0\,$.

On the other hand, for every $j\geq 1\,$,
\medskip

\noindent\hspace{1.13 cm}$\displaystyle \Big\|\,\frac 1{n_{j+1}-1}
\sum_{k=1}^{n_{j+1}-1}f_{\substack{{}\\ k}}\,
\Big\|_{\substack{{}\\2}}^{\, 2} =\frac 1{\big( n_{j+1}-1\big)^2}
\,\Big\|\, \sum_{l=1}^j\,\sum_{k=n_{\substack{{}\\ l}}}^{
n_{\substack{{}\\ {l+1}}}-1} f_{\substack{{}\\ k}}\,
\Big\|_{\substack{{}\\2}}^{\, 2}$

\noindent\hspace{4.773 cm}$\displaystyle
=\frac 1{\big( n_{j+1}-1\big)^2} \sum_{l=1}^j\big(
n_{\substack{{}\\ {l+1}}} - n_{\substack{{}\\ l}}\big)^2$
\medskip

\noindent\hspace{4.773 cm}$\displaystyle \geq \Big(\,
\frac{n_{j+1}-n_j}{n_{j+1}-1}\,\Big)^2 =
\Big(\, 1-\frac{n_j-1}{n_{j+1}-1}\,\Big)^2\geq \frac 14\; .$
\smallskip

\noindent Consequently, $\displaystyle \Big\|\,\frac 1n\,
\sum_{k=1}^n f_{\substack{{}\\ k}}\,\Big\|_2 \nrightarrow 0\,$,
so (\ref{unif.w.mix.}) does not hold for
$(f_{\substack{{}\\ k}})_{\substack{{}\\ k\geq 1}}\,$.

\end{proof}

In spite of the above examples, Theorem \ref{orbit} entails that
for orbits of power bounded linear operators weak mixing to zero
and uniformly weak mixing to zero are equivalent. We are now looking
for a larger class of vector sequences, for which weak mixing to
zero and uniformly weak mixing to zero are still equivalent.

Let us call a sequence
$(x_{\substack{{}\\ k}})_{\substack{{}\\ k\geq 1}}$ in a Banach
space $X$ \emph{convex shift-bounded} if there exists a constant $c>0$
such that
\begin{equation}\label{pow.bdd}
\Big\|\sum_{j=1}^p\,\lambda_j\, x_{j+k}\,\Big\|\leq c\,
\Big\|\sum_{j=1}^p\,\lambda_j\, x_j\,\Big\|\, ,\qquad k\geq 1
\end{equation}
holds for any choice of $p\in\mathbb{N}^*$ and
$\lambda_1\, ,\,\ldots\, ,\,\lambda_p\geq 0\,$. Clearly:
\smallskip

\noindent\hspace{0.5 cm}$\bullet\;\;$the convex shift-boundedness of a
sequence implies its boundedness;
\smallskip

\noindent\hspace{0.5 cm}$\bullet\;\;$if $U : X\longrightarrow X$ is a
power bounded linear operator and $x\in X\,$, then the

\noindent\hspace{0.9 cm}sequence
$\big( U^k(x)\big)_{\substack{{}\\ k\geq 1}}$ is convex shift-bounded.
\smallskip

We notice that not every convex shift-bounded sequence, even in
a Hilbert space, is the orbit of a bounded linear operator :

\begin{example}\label{ex3}
{\it Let us define the sequence
$(f_{\substack{{}\\ k}})_{\substack{{}\\ k\geq 1}}$ in
$L^2\big( [0,1]\big)$ by setting for every $k\in\mathbb N^*$ with
$k\equiv 1\, ({\rm mod}\, 4)$

\noindent\hspace{0.448 cm}$f_{\substack{{}\\ k}}(t):=t^{k}\, ,\quad
f_{\substack{{}\\ k+1}}(t):=t^{k+{\textstyle \frac 1{4\, (k+2)}}}
\, ,\quad
f_{\substack{{}\\ k+2}}(t):=t^{k+1}\, ,\quad
f_{\substack{{}\\ k+3}}(t):=t^{k+1+{\textstyle \frac 12}}\, .$
\medskip

\noindent Then $(f_{\substack{{}\\ k}})_{\substack{{}\\ k\geq 1}}$
is convex shift-bounded, but 
there exists no bounded linear operator $U : L^2\big( [0,1]\big)
\rightarrow L^2\big( [0,1]\big)$ such that
\begin{center}
$f_{\substack{{}\\ k}} =U^k(f)\, ,\qquad k\geq k_o$
\end{center}
for some $f\in L^2\big( [0,1]\big)$ and $k_o\in\mathbb N^*$.}
\end{example}

\begin{proof}
First of all, if $0<\alpha_{\substack{{}\\ 1}}<
\alpha_{\substack{{}\\ 2}}<\,\ldots\;$ are real numbers and
$g_{\substack{{}\\ k}}\in L^2\big( [0,1]\big)$ is defined by
$g_{\substack{{}\\ k}}(t):=
t^{\;\!\substack{{\alpha_{\substack{{}\\ k}}}\\{}}}$, then the
sequence $(g_{\substack{{}\\ k}})_{\substack{{}\\ k\geq 1}}$ is
convex shift-bounded. Indeed, for any $p\in\mathbb{N}^*$ and
$\lambda_1\, ,\,\ldots\, ,\,\lambda_p\geq 0\,$, the function
\medskip

\noindent\hspace{1.77 cm}$\displaystyle \mathbb N^*\ni k
\longmapsto \Big\|\sum_{j=1}^p\,\lambda_j\, g_{j+k}
\Big\|_2^{\, 2} =\sum_{j,j'=1}^p \lambda_j\,\lambda_{j'}\,
\frac 1{\alpha_{\substack{{}\\ j+k}} +
\alpha_{\substack{{}\\ j'+k}} +1}$
\smallskip

\noindent is decreasing. In particular, the sequence
$(f_{\substack{{}\\ k}})_{\substack{{}\\ k\geq 1}}$ is convex
shift-bounded.

On the other hand, if $\alpha\, ,\,\ep >0$ and we define
$h\in L^2\big( [0,1]\big)$ by $h(t):=t^\alpha -t^{\alpha +\ep}$,
then
\smallskip

\noindent\hspace{0.483 cm}$\displaystyle \| h\|_2^{\, 2} =
\int\limits_0^1\!\big(t^{2\alpha} +t^{2\alpha +2\ep} -
2\, t^{2\alpha +\ep}\big)\,{\rm d}t =
\frac{2\,\ep^2}{(2\alpha +1)\, (2\alpha +\ep +1)\,
(2\alpha +2\ep +1)}\; .$
\smallskip

\noindent It is easy to verify that
\begin{align}
&\| h\|_2^{\, 2}\leq\frac 1{2\, (2\alpha +1)}\,\Big(
\frac{\ep}{\alpha +\ep}\Big)^2\leq\frac 1{2\, (2\alpha +1)}\,
\Big(\frac{\ep}{\alpha}\Big)^2\,\hspace{0.61 cm}\text{ if }\,\ep
\leq 1
\; ,\label{upper} \\
&\| h\|_2^{\, 2}\geq\frac 1{4\, (2\alpha +1)}\,\Big(
\frac{\ep}{\alpha +\ep}\Big)^2\geq\frac 1{4\, (2\alpha +1)}\,
\Big(\frac{\ep}{\alpha +1}\Big)^2\,\text{ if }\,\ep\leq 1\, ,\,
\alpha\geq 2
\; .\label{lower}
\end{align}

Now let $k\in\mathbb N^*$ be arbitrary such that
$k\equiv 1\, ({\rm mod}\, 4)\,$. Then we have by (\ref{upper})
\medskip

\noindent\hspace{1.12 cm}$\displaystyle \| f_{\substack{{}\\ k}}
-f_{\substack{{}\\ k+1}}\|_2^{\, 2}\leq\frac 1{2\, (2k +1)}\,
\Big(\frac 1{4\, (k+2)\, k}\Big)^2 =
\frac 1{32\, k^2\, (k+2)^2\, (2k +1)}\, ,$
\medskip

\noindent while (\ref{lower}) yields
\medskip

\noindent\hspace{1.12 cm}$\displaystyle \| f_{\substack{{}\\ k+2}}
-f_{\substack{{}\\ k+3}} \|_2^{\, 2}\geq\frac 1{4\, (2k +3)}\,
\Big(\frac 1{2\, (k+2)}\Big)^2 =
\frac 1{16\, (k+2)^2\, (2k +3)}\, .$
\medskip

\noindent Consequently $\| f_{\substack{{}\\ k+2}}-
f_{\substack{{}\\ k+3}}\|_2^{\, 2}\geq k^2 \| f_{\substack{{}\\ k}}
-f_{\substack{{}\\ k+1}}\|_2^{\, 2}$, and so
\begin{equation}\label{antidomin.}
\| f_{\substack{{}\\ k+2}}-f_{\substack{{}\\ k+3}}
\|_{\substack{{}\\ 2}}\geq k\,\| f_{\substack{{}\\ k}}-
f_{\substack{{}\\ k+1}}\|_{\substack{{}\\ 2}}\, .
\end{equation}

Let us assume that there is a bounded linear operator
$U : L^2\big( [0,1]\big)\rightarrow L^2\big( [0,1]\big)$ such
that

\noindent\hspace{4.49 cm}$f_{\substack{{}\\ k}} =U^k(f)\, ,
\qquad k\geq k_o$
\medskip

\noindent for some $f\in L^2\big( [0,1]\big)$ and $k_o\in
\mathbb N^*$. Then, for every $k\geq k_o$ with $k\equiv 1\,
({\rm mod}\, 4)\,$, (\ref{antidomin.}) yields
\begin{center}
$k\,\| f_{\substack{{}\\ k}}-
f_{\substack{{}\\ k+1}}\|_{\substack{{}\\ 2}}\leq
\| f_{\substack{{}\\ k+2}}-f_{\substack{{}\\ k+3}}
\|_{\substack{{}\\ 2}} =\big\|\, U^2\big( f_{\substack{{}\\ k}}
-f_{\substack{{}\\ k+1}}\big)\big\|_{\substack{{}\\ 2}}
\leq\| U\|^2\,\| f_{\substack{{}\\ k}}-
f_{\substack{{}\\ k+1}}\|_{\substack{{}\\ 2}}\, ,$
\end{center}
hence $\,\| U\|\geq\sqrt{k}\,$. But this contradicts the
boundedness of $\, U\,$.

\end{proof}

We shall prove (in this section in the realm of reflexive
Banach spaces and in Section \ref{c.sh-bdd} in full generality)
that weak mixing to zero is equivalent with uniformly weak mixing
to zero for any convex shift-bounded sequence. First we prove an
easy implication of weak mixing to zero :

\begin{lemma}\label{zeroinclosure}
Let $(x_{\substack{{}\\ k}})_{\substack{{}\\ k\geq 1}}$ be a
bounded sequence in a Banach space $X\,$, which is weakly mixing
to zero, and $\mathcal A\subset\mathbb{N}^*$ with $D^*(\mathcal A
)>0\,$. Then the norm-closure of the convex hull
$\;{\rm conv} \big( \{ x_k\, ;\, k\in\mathcal A \}\big)$ of
$\,\{ x_k\, ;\, k\in\mathcal A \}$ contains $\, 0\,$.
\end{lemma}

\begin{proof}
Let us assume that $0$ is not in the norm-closure of
$\,{\rm conv} \big( \{ x_k\, ;\, k\in\mathcal A \}\big)\,$.
Then the Hahn-Banach theorem yields the existence of some
$\ep_o>0$ and $x^*\in X^*$ such that
\begin{equation}\label{separe}
\Re\, \langle\, x^*,x_{\substack{{}\\ k}}\rangle \geq \ep_o\, ,
\qquad k\in\mathcal A\, .
\end{equation}
Further, by a classical result of B. O. Koopman and J. von Neumann
(see e.g. \cite{Kr}, Chapter 2, (3.1) or \cite{NSZ}, Lemma 9.3),
there is a zero density set $\mathcal E\subset\mathbb{N}^*$ such
that
\begin{equation}\label{converge}
\lim_{\mathcal{E}\,\not\ni\, k\to\infty}\,
\langle\, x^*,x_{\substack{{}\\ k}}\rangle =0\, .
\end{equation}

Then $\mathcal A\,\backslash\,\mathcal E$ is infinite, because
otherwise we would get the contradiction
\smallskip

\noindent\hspace{3.407 cm}$0<D^*(\mathcal A)\leq D^*(\mathcal A
\,\backslash\,\mathcal E)+D^*(\mathcal E)=0\,$.
\smallskip

\noindent Let $\, k_{\substack{{}\\ 1}}<k_{\substack{{}\\ 2}} <\,
\ldots\,$ be the elements of $\mathcal A\,\backslash\,\mathcal E\,$.
Then (\ref{converge}) implies that $\,\langle\, x^*,
x_{\substack{{}\\ k_j}}\rangle\rightarrow 0\,$, in contradiction
with (\ref{separe}).

\end{proof}

For weakly relatively compact sequences a stronger statement holds,
which is essentially \cite{J2}, Corollary 2 :

\begin{lemma}\label{Jones}
A weakly relatively compact sequence
$(x_{\substack{{}\\ k}})_{\substack{{}\\ k\geq 1}}$ in a Banach
space $X$ is weakly mixing to zero if and only if there exists
a zero density set $\,\mathcal E\subset\mathbb{N}^*$ such that
\smallskip

\noindent\hspace{2.16 cm}$\displaystyle
\lim_{\mathcal{E}\,\not\ni\, k\to\infty} x_k =0\,$ with respect
to the weak topology of $X\,$.
\end{lemma}

\begin{proof}
An inspection of the proof of \cite{J2}, Corollary 2 shows that
it works for any weakly relatively compact sequence in a Banach
space, not only for those, which are orbits of power bounded
linear operators.

\end{proof}

We notice that, if
$(x_{\substack{{}\\ k}})_{\substack{{}\\ k\geq 1}}$ is a weakly
relatively compact sequence in a Banach space $X$, which is weakly
mixing to zero, and $\,\mathcal E\subset\mathbb{N}^*$ is as in
Lemma \ref{Jones}, then, according to the classical Mazur theorem
about the equality of the weak and norm closure of a convex subset
of $X$, the norm-closure of the convex hull of every infinite
subset of $\mathbb{N}^*\,\backslash\,\mathcal E$ contains $0\,$.
In particular, for any $\mathcal A\subset\mathbb{N}^*$ with
$D^*(\mathcal A)>0\,$, the norm-closure of the convex hull of the
infinite set $\mathcal A\,\backslash\,\mathcal E$ contains $0\,$.

Now we prove a consequence of the negation of uniformly weak mixing
to zero (cf. the first part of the proof of \cite{J1}, Theorem IV) :

\begin{lemma}\label{HB}
Let $(x_{\substack{{}\\ k}})_{\substack{{}\\ k\geq 1}}$ be a
sequence in the closed unit ball of a Banach space $X\,$, which is
not uniformly weakly mixing to zero. Then there exist
\begin{center}
$0<\ep_o\leq 1\,$,\\
$\mathcal B\subset\mathbb{N}^*$ with $\, D^*(\mathcal B)\geq
\ep_o\,$,\\
$k_1\, ,\, k_2\, ,\, \ldots\,\in\mathbb{N}^*$ with
$\, k_j-k_{j-1} > j\,$, \\
$x_1^*\, ,\, x_2^*\, ,\, \ldots\,\in X^*$ with $\,\| x_j^*\|\leq 1\,$,
\end{center}
such that
\begin{center}
$\displaystyle \mathcal B \cap\bigcup_{j\geq 2}
(\, k_{j-1}\, ,\, k_{j-1}+j\, ] =\emptyset\, ,$ \\
$\Re\,\langle\, x_j^*,x_{\substack{{}\\ k}}
\rangle > 2\,\ep_o\, ,\qquad k\in\mathcal B\cap
(\, k_{j-1}+j\, ,\, k_j\, ]\; ,\, j\geq 2\, .$
\end{center}
\end{lemma}

\begin{proof}
For any complex number $z$ we shall use the notations
\medskip

\noindent\hspace{1.6 cm}$\displaystyle \Re^+ z :=\begin{cases}
\Re\, z&\text{if $\,\Re\, z\geq 0$} \\
\; 0&\text{if $\,\Re\, z\leq 0$} \end{cases}\; ,\qquad
\Re^- z :=\begin{cases}
\; 0&\text{if $\,\Re\, z\geq 0$} \\
-\,\Re\, z&\text{if $\,\Re\, z\leq 0$} \end{cases}\; .$
\medskip

\noindent Then $\,\Re\, z=\Re^+ z -\Re^- z =\Re^+ z -\Re^+ (-z)\,$.

Since $(x_{\substack{{}\\ k}})_{\substack{{}\\ k\geq 1}}$ is not
uniformly weakly mixing to zero, there is $0<\ep_o\leq 1$ such that
\smallskip

\noindent\hspace{0.969 cm}$\displaystyle
\mathcal J := \bigg\{\, n\geq 1\, ;\, \sup\, \Big\{\, \frac 1n\,
\sum\limits_{k=1}^n\, |\langle\, x^*,x_{\substack{{}\\ k}}\rangle |\,
;\, x^*\in X^*\, ,\,\|x^*\|\leq 1\,\Big\} > 16\,\ep_o\,\bigg\}$
\smallskip

\noindent is infinite. Using (in the complex case)
$\langle x^*,x_{\substack{{}\\ k}}\rangle =\Re\,\langle\, x^*,
x_{\substack{{}\\ k}}\rangle -i\,\Re\,\langle\, i\, x^*,
x_{\substack{{}\\ k}}\rangle\, ,$
it follows that also
\smallskip

\noindent\hspace{0.789 cm}$\displaystyle
\mathcal J_{\substack{{}\\ \Re}} := \bigg\{\, n\geq 1\, ;\,
\sup\, \Big\{\,\frac 1n\,\sum\limits_{k=1}^n\, |\,\Re\,
\langle\, x^*,x_{\substack{{}\\ k}}\rangle |\, ;\, x^*\in X^*\, ,
\,\|x^*\|\leq 1\,\Big\} > 8\,\ep_o\,\bigg\}$
\smallskip

\noindent is infinite. Now, since $\,\Re\,\langle x^*,
x_{\substack{{}\\ k}}\rangle =\Re^+\langle\, x^*,
x_{\substack{{}\\ k}}\rangle -\Re^+\langle\, -\, x^*,
x_{\substack{{}\\ k}}\rangle\,$, we obtain that
\smallskip

\noindent\hspace{0.827 cm}$\displaystyle
\mathcal J_{\substack{{}\\ +}} := \bigg\{\, n\geq 1\, ;\,
\sup\, \Big\{\,\frac 1n\,\sum\limits_{k=1}^n\,\Re^+\langle\, x^*,
x_{\substack{{}\\ k}}\rangle\, ;\, x^*\in X^*\, ,\,\|x^*\|\leq 1\,
\Big\} > 4\,\ep_o\,\bigg\}$
\smallskip

\noindent is infinite.

Let $n\in\mathcal J_{\substack{{}\\ +}}$ be arbitrary. Then
there exists $y_n^*\in X^*$ with $\| y_n^*\|\leq 1$ such that
\smallskip

\noindent\hspace{4.36 cm}$\displaystyle
\frac 1n\,\sum\limits_{k=1}^n\,\Re^+\langle\, y_n^*,
x_{\substack{{}\\ k}}\rangle > 4\,\ep_o\, .$
\smallskip

\noindent Denoting $\,\mathcal B_n:=\{ 1\leq k\leq n\, ;\,
\Re^+\langle\, y_n^*,x_{\substack{{}\\ k}}\rangle >2\,\ep_o\}\,$,
we have
\medskip

\noindent\hspace{0.87 cm}$\displaystyle 4\,\ep_o <\frac 1n\,\Big(
\sum\limits_{k\in\mathcal B_n}\Re^+\langle\, y_n^*,
x_{\substack{{}\\ k}}\rangle +\sum_{\substack{1\leq k\leq n \\
k\notin\mathcal B_n}}\Re^+\langle\, y_n^*,
x_{\substack{{}\\ k}}\rangle\Big)\leq \frac 1n\;\mathrm{card}\,
(\mathcal B_n ) +2\,\ep_o\, ,$
\smallskip

\noindent hence $\,\mathrm{card}\, (\mathcal B_n )\geq 2\, n\,
\ep_o\,$.

Denoting now by $k_1$ the least element of
$\mathcal J_{\substack{{}\\ +}}\,$, we can construct recursively
a sequnce $k_1\, ,\, k_2\, ,\, \ldots\,\in
\mathcal J_{\substack{{}\\ +}}$ such that, for every $j\geq 2\,$,
\smallskip

\noindent\hspace{4.916 cm}$k_j-k_{j-1} > j\,$ and
\smallskip

\noindent\hspace{1.48 cm}the cardinality of
$\,{\mathcal B_{\substack{{}\\ k_j}}}' :=
\{\, k\in\mathcal B_{\substack{{}\\ k_j}}\, ;\, k>k_{j-1}+j\,\}\,$
is $\,\geq k_j\,\ep_o\,$.
\smallskip

\noindent Indeed, if we choose $k_j$ in the infinite set
$\mathcal J_{\substack{{}\\ +}}$ such that $k_j >
\frac{k_{j-1} +j}{\ep_o}\,$, then
\smallskip

\noindent\hspace{1.155 cm}$\mathrm{card}\, (
{\mathcal B_{\substack{{}\\ k_j}}}') \geq \mathrm{card}\, (
\mathcal B_{\substack{{}\\ k_j}}) - (k_{j-1} +j) \geq
2\, k_j\,\ep_o - (k_{j-1} +j) > k_j\,\ep_o\, .$
\medskip

Putting
\smallskip

\noindent\hspace{5.23 cm}$\displaystyle \mathcal B :=
\bigcup_{j\geq 2}\, {\mathcal B_{\substack{{}\\ k_j}}}'\, ,$
\smallskip

\noindent we have for every $j\geq 2$
\smallskip

\noindent\hspace{1.368 cm}$\mathcal B\cap (\, k_{j-1}\, ,\,
k_{j-1}+j\, ]=\emptyset\, ,\qquad
\mathcal B\cap (\, k_{j-1}+j\, ,\, k_j\, ] =
{\mathcal B_{\substack{{}\\ k_j}}}' \subset
\mathcal B_{\substack{{}\\ k_j}}\, ,$
\smallskip

\noindent and so

\noindent\hspace{1.393 cm}$\Re\,\langle\, y_{k_j}^{\, *},
x_{\substack{{}\\ k}}\rangle =\Re^+\langle\, y_{k_j}^{\, *},
x_{\substack{{}\\ k}}\rangle > 2\,\ep_o\,$ for all $\, k\in
\mathcal B\cap (\, k_{j-1}+j\, ,\, k_j\, ]\, .$
\medskip

\noindent On the other hand,
\medskip

\noindent\hspace{0.515 cm}$\displaystyle D^*(\mathcal B) =
\varlimsup_{n\to\infty}\,\frac 1n\;\mathrm{card}\, \big(
\mathcal B \cap [1,n]\big) \geq \varlimsup_{j\to\infty}\,
\frac 1{k_j}\;\mathrm{card}\, \big( \underbrace{\mathcal B \cap
(k_{j-1}+j,k_j]}_{=\,{\mathcal B_{\substack{{}\\ k_j}}}'}\big)
\geq\ep_o\, .$

\noindent Therefore, setting $x_j^{\;\! *}:=y_{k_j}^{\, *}\,$,
the proof is complete.

\end{proof}

We recall the following lemma of L. K. Jones on sequences of
integers (see \cite{J1}, Lemma 3 or \cite{J2}, Lemma) :

\begin{lemma}\label{weak-structure}
Let $\mathcal A_o\, ,\,\mathcal B$ be subsets of $\mathbb N^*$
with $D^*(\mathcal A_o)=1$ and $D^*(\mathcal B )>0\,$. Then
there exists an infinite subset $\mathcal I\subset
\mathcal A_o$ such that
\medskip

\noindent\hspace{2.31 cm}$D^*\big( \{ k\in\mathcal B\, ;\,
\mathcal F +k\subset\mathcal B\} \big) >0\,$ for any finite $\,
\mathcal F\subset\mathcal I\, .$
\end{lemma}
\hfill{$\square$}\par\medskip

Now, using the idea of the proof of \cite{J1}, Theorem IV,
we can prove that weak mixing to zero and uniformly weak
mixing to zero are equivalent for any convex shift-bounded
sequence in a reflexive Banach space :

\begin{proposition}\label{refl_wm=unif.wm}
For a convex shift-bounded sequence in a reflexive Banach space,
weak mixing to zero is equivalent to uniformly weak mixing to zero.
\end{proposition}

\begin{proof}
Let $(x_{\substack{{}\\ k}})_{\substack{{}\\ k\geq 1}}$ be a
convex shift-bounded sequence in the closed unit ball of a
reflexive Banach space $X\,$, which is weakly mixing to zero,
and let us assume that it is not uniformly weakly mixing to zero.
Let $c>0$ be such that (\ref{pow.bdd}) holds for any choice of
$p\in\mathbb{N}^*$ and $\lambda_1\, ,\,\ldots\, ,\,\lambda_p
\geq 0\,$.

By Lemma \ref{HB} there exist
\begin{center}
$0<\ep_o\leq 1\,$,\\
$\mathcal B\subset\mathbb{N}^*$ with $\, D^*(\mathcal B)\geq
\ep_o\,$,\\
$k_1\, ,\, k_2\, ,\, \ldots\,\in\mathbb{N}^*$ with
$\, k_j-k_{j-1} > j\,$, \\
$x_1^*\, ,\, x_2^*\, ,\,\ldots\,\in X^*$ with $\,\| x_j^*\|
\leq 1\,$,
\end{center}
such that

\noindent\hspace{4.064 cm}$\displaystyle \mathcal B \cap
\bigcup_{j\geq 2}(\, k_{j-1}\, ,\, k_{j-1}+j\, ]=\emptyset\, ,$
\smallskip

\noindent\hspace{2.2 cm}$\Re\,\langle\, x_j^*,
x_{\substack{{}\\ k}}\rangle > 2\,\ep_o\, ,\qquad k\in
\mathcal B\cap (\, k_{j-1}+j\, ,\, k_j\, ]\; ,\, j\geq 2\, .$
\medskip

On the other hand, since any bounded set in a reflexive Banach
space is weakly relatively compact, by Lemma \ref{Jones} there
exists $\mathcal A_o\subset\mathbb N^*$ with $D_*(\mathcal A_o)
=1$ such that $\displaystyle \lim_{\mathcal{A}_o\ni k\to\infty}
x_k =0\,$ with respect to the weak topology of $X\,$.

Finally, by Lemma \ref{weak-structure} there exists an infinite
subset $\mathcal I\subset\mathcal A_o$ such that
\medskip

\noindent\hspace{2.31 cm}$D^*\big( \{ k\in\mathcal B\, ;\,
\mathcal F +k\subset\mathcal B\} \big) >0\,$ for any finite $\,
\mathcal F\subset\mathcal I\, .$
\smallskip

Since $\displaystyle \lim_{\mathcal{I}\ni k\to\infty}
x_k =0\,$ with respect to the weak topology of $X$, there
are $p\in\mathbb{N}^*$, $n_{\substack{{}\\ 1}}<\,\ldots\, <n_p$
in $\mathcal I$ and $\lambda_{\substack{{}\\ 1}}\, ,\,\ldots\, ,
\,\lambda_p\geq 0\, ,\,\lambda_{\substack{{}\\ 1}}+ \ldots +
\lambda_p=1\,$, such that
\smallskip

\noindent\hspace{4.808 cm}$\displaystyle \Big\|\sum_{j=1}^p\,
\lambda_j\, x_{\substack{{}\\ n_j}}\Big\|\leq\frac{\ep_o}c\, .$
\medskip

\noindent By (\ref{pow.bdd}) it follows that
%
\begin{equation}\label{transl1}
\Big\|\sum_{j=1}^p\,\lambda_j\, x_{\substack{{}\\ {n_j +k}}}
\Big\|\leq c\,\Big\|\sum_{j=1}^p\,\lambda_j\,
x_{\substack{{}\\ n_j}}\Big\|\leq\ep_o\; ,\qquad k\geq 1\, .
\end{equation}

Now set $j_o:=\max\big( n_p-n_{\substack{{}\\ 1}},2\big)\,$.
Since the set $\big\{ k\in\mathcal B\, ;\,\{ n_{\substack{{}\\ 1}}
\, ,\,\ldots\, ,\, n_p\} +k\subset\mathcal B\big\}$ has strictly
positive upper density and so is infinite, it contains some $k$
such that $n_{\substack{{}\\ 1}} +k\geq k_{j_o}\,$. Then there
is a unique $j_{\substack{{}\\ 1}}\in\mathbb N^*$ with
$k_{j_{\substack{{}\\ 1}}-1}<n_{\substack{{}\\ 1}} +
k\leq k_{j_{\substack{{}\\ 1}}}\,$, for which we have $k_{j_o}
\leq n_{\substack{{}\\ 1}} +k\leq k_{j_{\substack{{}\\ 1}}}\,$,
hence $j_o\leq j_{\substack{{}\\ 1}}\,$. We claim that
\begin{equation}\label{lacunary1}
k_{j_{\substack{{}\\ 1}}-1}+j_{\substack{{}\\ 1}}<
n_{\substack{{}\\ 1}}+k\leq n_p+k\leq k_{j_{\substack{{}\\ 1}}}\; .
\end{equation}

Indeed, $k_{j_{\substack{{}\\ 1}}-1}<n_{\substack{{}\\ 1}} +k\,$,
$n_{\substack{{}\\ 1}} +k\in\mathcal B\,$ and $\,\mathcal B\cap
(k_{j_{\substack{{}\\ 1}}-1},k_{j_{\substack{{}\\ 1}}-1}+
j_{\substack{{}\\ 1}}]=\emptyset\,$ imply that
$k_{j_{\substack{{}\\ 1}}-1}+j_{\substack{{}\\ 1}}<
n_{\substack{{}\\ 1}} +k\,$.
Similarly, $n_p+k =n_{\substack{{}\\ 1}} +k +(n_p -
n_{\substack{{}\\ 1}} )\leq k_{j_{\substack{{}\\ 1}}} +j_o
<k_{j_{\substack{{}\\ 1}}} +j_{\substack{{}\\ 1}} +1\,$,
$n_p+k\in\mathcal B\,$ and $\,\mathcal B\cap
(k_{j_{\substack{{}\\ 1}}}\! ,k_{j_{\substack{{}\\ 1}}} +
j_{\substack{{}\\ 1}} +1]=\emptyset\,$ yield $\, n_p+k\leq
k_{j_{\substack{{}\\ 1}}}\,$.

By (\ref{lacunary1}) we have $\, n_j+k\in\mathcal B\cap
(k_{j_{\substack{{}\\ 1}}-1}+j_{\substack{{}\\ 1}},
k_{j_{\substack{{}\\ 1}}}]\, , 1\leq j\leq p\,$, so
\medskip

\noindent\hspace{3.456 cm}$\Re\,\langle\,
x_{j_{\substack{{}\\ 1}}}^{\;\! *}, x_{\substack{{}\\ n_j+k}}
\rangle > 2\,\ep_o\, ,\qquad 1\leq j\leq p\, .$
\smallskip

\noindent Since $\| x_{j_{\substack{{}\\ 1}}}^{\;\! *}\|\leq 1\,$,
it follows that
\smallskip

\noindent\hspace{0.847 cm}$\displaystyle \Big\|\sum_{j=1}^p\,
\lambda_j\, x_{\substack{{}\\ {n_j +k}}}\Big\|\geq \Re\,\left<
x_{j_{\substack{{}\\ 1}}}^{\;\! *}\, ,\,\sum_{j=1}^p\,\lambda_j\,
x_{\substack{{}\\ {n_j +k}}} \right> =\sum_{j=1}^p\,\lambda_j\,
\Re\,\langle\, x_{j_{\substack{{}\\ 1}}}^{\;\! *},
x_{\substack{{}\\ k_j}}\rangle > 2\,\ep_o\, ,$
\smallskip

\noindent in contradiction with (\ref{transl1}).

\end{proof}

If in Lemma \ref{weak-structure} the set $\mathcal I$ would be
not only infinite, but with $D^*(\mathcal I)>0\,$, then in the
proof of Proposition \ref{refl_wm=unif.wm} we could use Lemma
\ref{zeroinclosure} instead of Lemma \ref{Jones} and so we
would get a proof of Proposition \ref{refl_wm=unif.wm} without
the reflexivity assumption. In the next section we shall prove
a result like Lemma \ref{weak-structure} (Theorem \ref{structure}),
which implies that for
every $\mathcal B\subset\mathbb N^*$ with $D^*(\mathcal B)>0$
there is a set $\mathcal A\subset\mathbb N^*$ with
$D^*(\mathcal A)>0$ such that any finite subset of $\mathcal A$
has infinitely many translates contained in $\mathcal B\,$. This
result will enable us to eliminate the reflexivity condition in
Proposition \ref{refl_wm=unif.wm}.

\bigskip
\section{Sets of non-zero upper Banach density}

We recall that the \emph{upper Banach density} $BD^*(\mathcal B)$
of some $\mathcal B\subset\mathbb N=N^*\cup\{ 0\}$ is defined by
\begin{equation*}
BD^*(\mathcal B):=
\varlimsup_{\substack{a,b\in\mathbb N\\ b-a\to\infty}}\,
\frac 1{b-a+1}\;\mathrm{card}\, \big(\mathcal B \cap [a,b]\big)
=\varlimsup_{\substack{a,b\in\mathbb N^*\\ b-a\to\infty}}\,
\frac 1{b-a+1}\;\mathrm{card}\, \big(\mathcal B \cap [a,b]\big)
\end{equation*}
(see e.g. \cite{F2}, Chapter 3, \S 5). For any $\mathcal B\subset
\mathbb N^*$ we have $BD^*(\mathcal B)\geq D^*(\mathcal B)\,$,
but it is easily seen that $BD^*(\mathcal B) > D^*(\mathcal B)$
can happen. In this section we investigate the structure of the
sets $\mathcal B\subset\mathbb N^*$ with $BD^*(\mathcal B)>0$ by
proving a precise version of the theorem of R. Ellis \cite{F2},
Theorem 3.20. The proof is based on the ergodic theoretical
methods of H. Furstenberg exposed in \cite{F2}, Chapter 3, \S 5.

Let us consider $\Omega :=\{ 0,1\}^{\mathbb N}$ and endow it with
the metrizable compact product topology of the discrete topologies
on $\{ 0,1\}\,$. We shall denote the components of $\omega\in
\Omega$ by $\omega_{\substack{{}\\ k}}\,$, so that $\omega =
(\omega_{\substack{{}\\ k}})_{\substack{{}\\ k\in\mathbb N}}\,$.
For every $\mathcal B\subset\mathbb N$ we define
$\omega^{(\mathcal B)}\in\Omega$ by setting
\begin{center}
$\omega^{(\mathcal B)}_k:=\begin{cases}
\; 1 &\text{if $\, k\in\mathcal B$}\\
\; 0 &\text{if $\, k\notin\mathcal B$} \end{cases}\; .$
\end{center}
In other words, $\,\omega^{(\mathcal B)}$ is the characteristic
function of $\mathcal B\,$, considered an element of $\Omega\,$.
Clearly, $\,\mathcal B\longmapsto\omega^{(\mathcal B)}$ is a
bijection of the set of all subsets of $\mathbb N$ onto $\Omega\,$.

Let $s_{\substack{{}\\ \!\leftarrow}}$ denote the backward
shift on $\Omega\,$, defined by
\smallskip

\noindent\hspace{4.263 cm}$s_{\substack{{}\\ \!\leftarrow}}
\big( (\omega_{\substack{{}\\ k}})_{\substack{{}\\ k\in\mathbb N}}
\big) =(\omega_{\substack{{}\\ k+1}})_{\substack{{}\\ k\in
\mathbb N}}\,$,
\smallskip

\noindent and set, for every $\mathcal B\subset\mathbb N\,$,
\smallskip

\noindent\hspace{4.188 cm}$\Omega^{(\mathcal B)}:=\overline{
\{ s_{\substack{{}\\ \!\leftarrow}}^{\, n}(\omega^{(\mathcal B)})
\, ;\, n\geq 0\}}\;$.
\smallskip

\noindent Clearly, $s_{\substack{{}\\ \!\leftarrow}}
(\Omega^{(\mathcal B)})\subset\Omega^{(\mathcal B)}\,$.

The following result is the one-sided version of \cite{F2},
Lemma 3.17 and it establishes a link between upper Banach
density and the ergodic theory of the dynamical system
$(\Omega , s_{\substack{{}\\ \!\leftarrow}})\,$. Its proof is
almost identical to the proof of \cite{F2}, Lemma 3.17 and we
sketch it only for the sake of completeness :

\begin{lemma}\label{erg.inv.meas}
For every $\mathcal B\subset\mathbb N$ and every $\ep >0$
there exists an ergodic
$s_{\substack{{}\\ \!\leftarrow}}$-invariant probability Borel
measure $\mu$ on $\Omega^{(\mathcal B)}$ such that
\begin{center}
$\mu\big( \{\omega\in\Omega^{(\mathcal B)}\, ;\,
\omega_{\substack{{}\\ 0}} =1\}\big) > BD^*(\mathcal B) -\ep
\, .$
\end{center}
\end{lemma}

\begin{proof}
Choose some $\, a_{\substack{{}\\ 1}}\, ,\,
b_{\substack{{}\\ 1}}\, ,\, a_{\substack{{}\\ 2}}\, ,\,
b_{\substack{{}\\ 2}}\, ,\;\ldots\;\in\mathbb N\,$ with $\,
b_j-a_j\geq j\, ,\, j\geq 0\, ,$ such that
%
%
\smallskip

\noindent\hspace{3 cm}$\displaystyle \frac 1{b_j-a_j+1}\,
\mathrm{card}\, \big(\mathcal B \cap [a_j,b_j]\big)
\longrightarrow BD^*(\mathcal B)\, .$
\medskip

\noindent Passing to a subsequence if necessary, we can assume
that, for any continuous function $f\in C\big(\Omega^{(\mathcal B)}
\big)\,$,

\noindent\hspace{1.78 cm}the limit $\displaystyle I(f):=
\lim_{j\to\infty}\,\frac 1{b_j-a_j+1}\,\sum_{n=a_j}^{b_j}
f\big( s_{\substack{{}\\ \!\leftarrow}}^{\, n}
(\omega^{(\mathcal B)})\big)\,$ exists.
\smallskip

\noindent Then $I$ is a positive linear functional on $C\big(
\Omega^{(\mathcal B)}\big)$ and $I(1)=1\,$. Moreover,
\begin{equation}\label{invariance}
I(f\circ s_{\substack{{}\\ \!\leftarrow}})=I(f)\, ,\qquad
f\in C\big(\Omega^{(\mathcal B)}\big)\, .
\end{equation}
Indeed, for every $f\in C\big(\Omega^{(\mathcal B)}\big)\,$,

\noindent\hspace{0.23 cm}$\displaystyle I(f\circ
s_{\substack{{}\\ \!\leftarrow}}) -I(f)
=\lim_{j\to\infty}\,\frac 1{b_j-a_j+1}\,\Big(\,\sum_{n=a_j}^{b_j}
f\big( s_{\substack{{}\\ \!\leftarrow}}^{\, n+1}
(\omega^{(\mathcal B)})\big) -\sum_{n=a_j}^{b_j}
f\big( s_{\substack{{}\\ \!\leftarrow}}^{\, n}
(\omega^{(\mathcal B)})\big) \Big)$
\smallskip

\noindent\hspace{2.834 cm}$\displaystyle =\lim_{j\to\infty}\,
\frac 1{b_j-a_j+1}\,\Big(\, f\big(
s_{\substack{{}\\ \!\leftarrow}}^{\, b_j+1}(\omega^{(\mathcal B)})
\big) -f\big( s_{\substack{{}\\ \!\leftarrow}}^{\, a_j}
(\omega^{(\mathcal B)})\big) \Big)$
\smallskip

\noindent\hspace{2.834 cm}$\displaystyle \leq\lim_{j\to\infty}\,
\frac{2\,\| f\|_\infty}{b_j-a_j+1} =0\, .$
\medskip

\noindent By the Riesz representation theorem there exists a
probability Borel measure $\nu_{\substack{{}\\ I}}$ on
$\Omega^{(\mathcal B)}$ such that
\smallskip

\noindent\hspace{3.153 cm}$\displaystyle I(f)=
\int\limits_{\Omega^{(\mathcal B)}} f(\omega )\,\mathrm{d}
\nu_{\substack{{}\\ I}}(\omega)\, ,\qquad
f\in C\big(\Omega^{(\mathcal B)}\big)\, .$
\smallskip

\noindent Property (\ref{invariance}) of $I$ implies that
$\nu_{\substack{{}\\ I}}$ is
$s_{\substack{{}\\ \!\leftarrow}}$-invariant. Moreover, since the
characteristic function $\chi$ of $\{\omega\in\Omega^{(\mathcal B)}
\, ;\,\omega_{\substack{{}\\ 0}} =1\}$ is continuous, we have

\noindent\hspace{0.2 cm}$\displaystyle \nu_{\substack{{}\\ I}}\big(
\{\omega\in\Omega^{(\mathcal B)}\, ;\,\omega_{\substack{{}\\ 0}} =1
\}\big) =\int\limits_{\Omega^{(\mathcal B)}}\chi (\omega )\,
\mathrm{d}\nu_{\substack{{}\\ I}}(\omega) =\lim_{j\to\infty}\,
\frac 1{b_j-a_j+1}\,\sum_{n=a_j}^{b_j}\chi\big(
s_{\substack{{}\\ \!\leftarrow}}^{\, n}(\omega^{(\mathcal B)})\big)$
\smallskip

\noindent\hspace{3.95 cm}$\displaystyle =\lim_{j\to\infty}\,
\frac 1{b_j-a_j+1}\,\mathrm{card}\,\big(\mathcal B \cap [a_j,b_j]
\big) =BD^*(\mathcal B)\, .$
\medskip

The convex set $\mathcal P^{s_{\substack{{}\\ \!\leftarrow}}}\big(
\Omega^{(\mathcal B)}\big)$ of all
$s_{\substack{{}\\ \!\leftarrow}}$-invariant probability Borel
measures on $\Omega^{(\mathcal B)}\,$, considered imbedded in
the dual space of $C\big(\Omega^{(\mathcal B)}\big)\,$, is
weak${}^*$-compact and its extreme points are the ergodic measures
in $\mathcal P^{s_{\substack{{}\\ \!\leftarrow}}}\big(
\Omega^{(\mathcal B)}\big)$
(see, for example, \cite{F2}, Proposition 3.4). According to the
Krein-Milman theorem, it follows that $\nu_{\substack{{}\\ I}}$
is a weak${}^*$-limit of convex combinations of ergodic meaures in
$\mathcal P^{s_{\substack{{}\\ \!\leftarrow}}}\big(
\Omega^{(\mathcal B)}\big)\,$. Therefore, since
$\nu_{\substack{{}\\ I}}\big(\{\omega\in\Omega^{(\mathcal B)}\, ;
\,\omega_{\substack{{}\\ 0}} =1\}\big) =BD^*(\mathcal B)\,$, we
conclude that there exists an ergodic measure $\mu\in
\mathcal P^{s_{\substack{{}\\ \!\leftarrow}}}\big(
\Omega^{(\mathcal B)}\big)$ such that $\mu\big(\{\omega\in
\Omega^{(\mathcal B)}\, ;\,\omega_{\substack{{}\\ 0}} =1\}\big)
>BD^*(\mathcal B)-\ep\,$.

\end{proof}

Now we prove the announced extension of \cite{F2}, Theorem 3.20 :

\begin{theorem}\label{structure}
If $\mathcal B\subset\mathbb N$ and $0<\ep <BD^*(\mathcal B)\,$,
then there exist
\smallskip
\begin{center}
$\mathcal A\subset\mathbb N\,$ having density $\, D(\mathcal A)
>BD^*(\mathcal B)-\ep\, ,$ \\
$0\leq m_{\substack{{}\\ 1}}<m_{\substack{{}\\ 2}}<\,\ldots\;$
and $\, 0\leq n_{\substack{{}\\ 1}}<n_{\substack{{}\\ 2}}<
\,\ldots\;$ in $\,\mathbb N\, ,$
\end{center}
for which
\begin{center}
$\mathcal A\cap [0,m_j] =\big\{ k\in [0,m_j]\, ;\, k+n_j\in
\mathcal B\big\}\; ,\qquad j\geq 1\, ,$
\end{center}
that is
\begin{center}
$k\in\mathcal A\,\Longleftrightarrow\, k+n_j\in\mathcal B\,$
whenever $\; 0\leq k\leq m_j\, ,\, j\geq 1\, .$
\end{center}
\end{theorem}

\begin{proof}
For every $\omega\in\Omega$ we set $\mathcal A_\omega
=\{ k\in\mathbb N\, ;\,\omega_k=1\}\,$, so that $\omega =
\omega^{(\mathcal A_\omega )}\,$. Clearly,
$\mathcal A_{\omega^{(\mathcal B)}}=\mathcal B\,$.

By Lemma \ref{erg.inv.meas} there exists an ergodic
$s_{\substack{{}\\ \!\leftarrow}}$-invariant probability Borel
measure $\mu$ on $\Omega^{(\mathcal B)}$ such that
\smallskip

\noindent\hspace{2.832 cm}$\mu_{\substack{{}\\ \mathcal B}}:=
\mu\big(\{\omega\in\Omega^{(\mathcal B)}\, ;\,
\omega_{\substack{{}\\ 0}} =1\}\big) > BD^*(\mathcal B) -\ep\, .$
\medskip

\noindent Let $\chi$ denote the characteristic function of
$\{\omega\in\Omega^{(\mathcal B)}\, ;\,\omega_{\substack{{}\\ 0}}
=1\}\subset\Omega^{(\mathcal B)}\,$. Then, by
the Birkhoff ergodic theorem, for $\mu$-almost every $\omega\in
\Omega^{(\mathcal B)}$ we have
\begin{equation}\label{Birkh.}
\frac 1{n+1}\;\mathrm{card}\,\big(\mathcal A_\omega \cap [0,n]\big)
=\frac 1{n+1}\,\sum_{k=0}^{n}\,\chi\big(
s_{\substack{{}\\ \!\leftarrow}}^{\, k}(\omega)\big)
\longrightarrow\mu_{\substack{{}\\ \mathcal B}}\, .
\end{equation}
%
%
\noindent Let $\Omega^{(\mathcal B)}_{\rm Birkhoff}$ be the set
of all $\omega\in\Omega^{(\mathcal B)}$, for which (\ref{Birkh.})
holds. Then
\smallskip

\noindent\hspace{1 cm}$\mathcal A_\omega$ has density
$D(\mathcal A_\omega )=\mu_{\substack{{}\\ \mathcal B}}
> BD^*(\mathcal B) -\ep\,$ for every $\,\omega\in
\Omega^{(\mathcal B)}_{\rm Birkhoff}\;$,
\smallskip
%

\noindent\hspace{1 cm}$\Omega^{(\mathcal B)}_{\rm Birkhoff}\,$
is $\mu$-measurable and $\,\mu\big(\Omega^{(\mathcal B)}\,
\backslash\,\Omega^{(\mathcal B)}_{\rm Birkhoff}\big) =0\;$.
\medskip

\underline{Case 1}: there exists $\omega\in
\Omega^{(\mathcal B)}_{\rm Birkhoff}\,\backslash\,
\{ s_{\substack{{}\\ \!\leftarrow}}^{\, n}
(\omega^{(\mathcal B)})\, ;\, n\geq 0\}\,$.
\smallskip

Set $\mathcal A:=\mathcal A_\omega$ and choose some
$m_{\substack{{}\\ 1}}\geq 0\,$. Since
\begin{equation}\label{appr.}
\omega\in\Omega^{(\mathcal B)}=\overline{
\{ s_{\substack{{}\\ \!\leftarrow}}^{\, n}(\omega^{(\mathcal B)})
\, ;\, n\geq 0\}}\, ,
\end{equation}
there exists a smallest $n_{\substack{{}\\ 1}}\geq 0$ such that
\smallskip

\noindent\hspace{1.175 cm}$\omega_{\substack{{}\\ k}}=
s_{\substack{{}\\ \!\leftarrow}}^{\, n_1}
(\omega^{(\mathcal B)})_{\substack{{}\\ k}}=
\omega^{(\mathcal B)}_{\substack{{}\\ k+n_1}}=
\begin{cases}
\; 1 &\text{if $\, k+n_{\substack{{}\\ 1}}\in\mathcal B$}\\
\; 0 &\text{if $\, k+n_{\substack{{}\\ 1}}\notin\mathcal B$}
\end{cases}\; ,\qquad 0\leq k\leq m_{\substack{{}\\ 1}}\, ,$
\smallskip

\noindent that is $\mathcal A\cap [0,m_{\substack{{}\\ 1}}] =
\big\{ k\in [0,m_{\substack{{}\\ 1}}]\, ;\,
k+n_{\substack{{}\\ 1}}\in\mathcal B\big\}\,$.

Next $\omega\neq s_{\substack{{}\\ \!\leftarrow}}^{\, n_1}
(\omega^{(\mathcal B)})$ implies that
$\omega_{\substack{{}\\ m_{\substack{{}\\ 2}}}}\neq
s_{\substack{{}\\ \!\leftarrow}}^{\, n_1}
(\omega^{(\mathcal B)})_{\substack{{}\\ m_{\substack{{}\\ 2}}}}$
for some $m_{\substack{{}\\ 2}}\in\mathbb N\,$. Since
$\omega_{\substack{{}\\ k}}=
s_{\substack{{}\\ \!\leftarrow}}^{\, n_1}
(\omega^{(\mathcal B)})_{\substack{{}\\ k}}$ for all $0\leq k
\leq m_{\substack{{}\\ 1}}\,$, we have $m_{\substack{{}\\ 1}}
<m_{\substack{{}\\ 2}}\,$. Now, again by (\ref{appr.}), there
exists a smallest $n_{\substack{{}\\ 2}}\geq 0$ such that
\smallskip

\noindent\hspace{1.175 cm}$\omega_{\substack{{}\\ k}}=
s_{\substack{{}\\ \!\leftarrow}}^{\, n_2}
(\omega^{(\mathcal B)})_{\substack{{}\\ k}}=
\omega^{(\mathcal B)}_{\substack{{}\\ k+n_2}}=
\begin{cases}
\; 1 &\text{if $\, k+n_{\substack{{}\\ 2}}\in\mathcal B$}\\
\; 0 &\text{if $\, k+n_{\substack{{}\\ 2}}\notin\mathcal B$}
\end{cases}\; ,\qquad 0\leq k\leq m_{\substack{{}\\ 2}}\, ,$
\smallskip

\noindent that is $\mathcal A\cap [0,m_{\substack{{}\\ 2}}] =
\big\{ k\in [0,m_{\substack{{}\\ 2}}]\, ;\,
k+n_{\substack{{}\\ 2}}\in\mathcal B\big\}\,$. By the minimality
property of $n_{\substack{{}\\ 1}}$ we have $n_{\substack{{}\\ 1}}
\leq n_{\substack{{}\\ 2}}\,$, while
$\,\omega_{\substack{{}\\ m_{\substack{{}\\ 2}}}}\neq
s_{\substack{{}\\ \!\leftarrow}}^{\, n_1}
(\omega^{(\mathcal B)})_{\substack{{}\\ m_{\substack{{}\\ 2}}}}$ 
yields $n_{\substack{{}\\ 1}}\neq n_{\substack{{}\\ 2}}\,$.
Therefore $n_{\substack{{}\\ 1}} < n_{\substack{{}\\ 2}}\,$.

By induction we obtain $m_{\substack{{}\\ 1}}<m_{\substack{{}\\ 2}}
<\,\ldots\;$ and $\, n_{\substack{{}\\ 1}}<n_{\substack{{}\\ 2}}<
\,\ldots\;$ in $\mathbb N$ such that
\smallskip

\noindent\hspace{2.18 cm}$\mathcal A\cap [0,m_j] =\big\{ k\in
[0,m_j]\, ;\, k+n_j\in\mathcal B\big\}\,$ for all $\, j\geq 1\,$.
\medskip

\underline{Case 2}: $\Omega^{(\mathcal B)}_{\rm Birkhoff}\subset
\{ s_{\substack{{}\\ \!\leftarrow}}^{\, n}
(\omega^{(\mathcal B)})\, ;\, n\geq 0\}\,$.
\smallskip

We claim that there exists a smallest $n_o\in\mathbb N^*$ such that
$s_{\substack{{}\\ \!\leftarrow}}^{\, n_o}(\omega^{(\mathcal B)})
=\omega^{(\mathcal B)}\,$.

For let us assume that all $s_{\substack{{}\\ \!\leftarrow}}^{\, n}
(\omega^{(\mathcal B)})$ are different. Then, for every $n\geq 0\,$,
since
\smallskip

\noindent\hspace{1.29 cm}$\{ s_{\substack{{}\\ \!\leftarrow}}^{\, n}
(\omega^{(\mathcal B)})\}\subset
s_{\substack{{}\\ \!\leftarrow}}^{\, -1}\big(
s_{\substack{{}\\ \!\leftarrow}}^{\, n+1}(\omega^{(\mathcal B)})\big)
\subset\{ s_{\substack{{}\\ \!\leftarrow}}^{\, n}
(\omega^{(\mathcal B)})\}\cup \big(\Omega^{(\mathcal B)}\,
\backslash\,\Omega^{(\mathcal B)}_{\rm Birkhoff}\big)\, ,$
\smallskip

\noindent by the $s_{\substack{{}\\ \!\leftarrow}}$-invariance
of $\mu$ we obtain $\mu\big( \{ s_{\substack{{}\\ \!\leftarrow}}^{\, n}
(\omega^{(\mathcal B)})\}\big) =\mu\big( \{
s_{\substack{{}\\ \!\leftarrow}}^{\, n+1}(\omega^{(\mathcal B)})\}
\big)\,$. Thus
\smallskip

\noindent\hspace{1.723 cm}$\displaystyle \mu (\Omega^{(\mathcal B)})
=\sum_{n=0}^\infty\,\mu\big( \{ s_{\substack{{}\\ \!\leftarrow}}^{\, n}
(\omega^{(\mathcal B)})\}\big) =\begin{cases}
\;\; 0&\text{if $\,\mu\big( \{\omega^{(\mathcal B)})\}\big) =0$} \\
+\infty&\text{if $\,\mu\big( \{\omega^{(\mathcal B)})\}\big) >0$}
\end{cases}\, ,$
\smallskip

\noindent in contradiction with $\mu (\Omega^{(\mathcal B)})=1\,$.

Now $s_{\substack{{}\\ \!\leftarrow}}^{\, n_o}(\omega^{(\mathcal B)})
=\omega^{(\mathcal B)}$ means that $k\in\mathbb N$ belongs to
$\mathcal B$ if and only if $k+n_o\in\mathcal B\,$. Therefore,
with $\mathcal A:=\mathcal B\,$, any $0\leq m_{\substack{{}\\ 1}}
<m_{\substack{{}\\ 2}} <\,\ldots\;$ and $n_j:=j\, n_o\,$, we have
\smallskip

\noindent\hspace{2.86 cm}$\displaystyle D(\mathcal A)=
\frac 1{n_o}\;\mathrm{card}\,\big(\mathcal B \cap [0,n_o -1]\big)
=BD^*(\mathcal B)\,$,
\smallskip

\noindent\hspace{2.333 cm}$\mathcal A\cap [0,m_j] =
\big\{ k\in [0,m_j]\, ;\, k+n_j\in
\mathcal B\big\}\; ,\qquad j\geq 1\, .$

\end{proof}

We recall that a celebrated theorem of E. Szemer\'edi (answering
a conjecture of P. Erd\H{o}s) states that if $\mathcal B\subset
\mathbb N^*$ has non-zero upper Banach density, then it
contains arbitrarily long arithmetic progressions \cite{Sz}.
H. Furstenberg gave a new ergodic theoretical proof of
Szemer\'edi's theorem by deducing it from a far-reaching
multiple recurrence theorem \cite{F1} (see also \cite{F2},
Chapter 3, \S 7). It is interesting to notice, even if it
looks not to be relevant, that via Theorem \ref{structure} the
proof of Szemer\'edi's theorem can be reduced to the case when
$\mathcal B$ has non-zero density.

The above theorem implies the following counterpart of Lemma
\ref{weak-structure} :

\begin{corollary}\label{suff.-structure}
Let $\mathcal A_o\, ,\,\mathcal B$ be subsets of $\mathbb N^*$
with $D^*(\mathcal A_o)=1$ and $0<\ep <BD^*(\mathcal B )\,$.
Then there exists $\mathcal I\subset \mathcal A_o$ with
$D^*(\mathcal I)>BD^*(\mathcal B ) -\ep\,$, such that
\medskip

\noindent\hspace{2.214 cm}$\{ k\in\mathbb N\, ;\,\mathcal F +k
\subset\mathcal B\}\,$ is infinite for any finite $\,
\mathcal F\subset\mathcal I\, .$
\end{corollary}

\begin{proof}
By Theorem \ref{structure} there exist $\mathcal A\subset
\mathbb N\,$ having density $\, D(\mathcal A) >BD^*(\mathcal B)
-\ep\,$, as well as
$\, 0\leq m_{\substack{{}\\ 1}}<m_{\substack{{}\\ 2}}<\,\ldots\;$
and $\, 0\leq n_{\substack{{}\\ 1}}<n_{\substack{{}\\ 2}}<
\,\ldots\;$ in $\,\mathbb N\,$, such that
\smallskip

\noindent\hspace{2.333 cm}$\mathcal A\cap [0,m_j] =
\big\{ k\in [0,m_j]\, ;\, k+n_j\in
\mathcal B\big\}\; ,\qquad j\geq 1\, .$
\smallskip

Set $\mathcal I :=\mathcal A\cap \mathcal A_o\,$. Since
\smallskip

\noindent\hspace{1.53 cm}$1=D^*(\mathcal A_o)\leq
D^*(\mathcal A\cap \mathcal A_o)+D^*(\mathbb N\,\backslash\,
\mathcal A) =D^*(\mathcal I) +1-D(\mathcal A)\, ,$
\smallskip

\noindent we have $D^*(\mathcal I)\geq D(\mathcal A)>BD^*(
\mathcal B)-\ep\,$. On the other hand, for any $j\geq 1\,$,
the set
\smallskip

\noindent\hspace{3.922 cm}$\big\{ k\in\mathbb N\, ;\,\big(
\mathcal I\cap [0,m_j]\big)+k\subset\mathcal B\big\}$
\smallskip

\noindent contains $\{ n_j\, ,\, n_{j+1}\, ,\,\ldots\,\}\,$,
hence is infinite.

\end{proof}

\bigskip
\section{Weak mixing to zero for convex shift-bounded sequences}
\label{c.sh-bdd}

Using Theorem \ref{structure}, in this section we show that
Proposition \ref{refl_wm=unif.wm} holds without the reflexivity
assumption. Actually we shall prove a slightly more general
result, stating that any convex shift-bounded sequence in a
Banach space, which is weakly mixing to zero, satisfies
(\ref{B.unif.w.mix.}). For the proof we shall use the following
counterpart of Lemma \ref{HB} for the sequences not satisfying
(\ref{B.unif.w.mix.}) :

\begin{lemma}\label{BHB}
Let $(x_{\substack{{}\\ k}})_{\substack{{}\\ k\geq 1}}$ be a
sequence in the closed unit ball of a Banach space $X\,$, such
that, for some
$a_{\substack{{}\\ 1}}\, ,\, b_{\substack{{}\\ 1}}\, ,\,
a_{\substack{{}\\ 2}}\, ,\, b_{\substack{{}\\ 2}}\, ,\,\ldots\,
\in\mathbb N^*\,$ with $\, b_j-a_j\geq j\, ,\, j\geq 1\, ,$
we have
\begin{center}
$\displaystyle \varlimsup_{j\to\infty}\;
\sup \Big\{\, \frac 1{b_j-a_j+1}\, \sum\limits_{k=a_j}^{b_j}\,
|\langle x^*, x_{\substack{{}\\ k}}\rangle |\, ;\, x^*\in X^*
\, ,\,\|x^*\|\leq 1\,\Big\} >0\, .$
\end{center}
Then there exist
\begin{center}
$0<\ep_o\leq 1\,$,\\
$\mathcal B\subset\mathbb{N}^*$ with $\, BD^*(\mathcal B)\geq
\ep_o\,$,\\
$j_{\substack{{}\\ 1}} < j_{\substack{{}\\ 2}} <\, \ldots\,$
in $\,\mathbb{N}^*\!$ with $\, b_{j_n}-b_{j_{n-1}} > n\,$, \\
$x_1^*\, ,\, x_2^*\, ,\, \ldots\,\in X^*$ with $\,\| x_n^*\|
\leq 1\,$,
\end{center}
such that
\begin{center}
$\displaystyle \mathcal B \cap\bigcup_{n\geq 2}
(\, b_{j_{n-1}}\, ,\, b_{j_{n-1}}+n\, ] =\emptyset\, ,$ \\
$\Re\,\langle\, x_n^*,x_{\substack{{}\\ k}}
\rangle > 2\,\ep_o\, ,\qquad k\in\mathcal B\cap
(\, b_{j_{n-1}}+n\, ,\, b_{j_n}\, ]\; ,\, n\geq 2\, .$
\end{center}
\end{lemma}

\begin{proof}
We shall proceed similarly as in the proof of Lemma \ref{HB}.

Let $0<\ep_o\leq 1$ be such that

\noindent\hspace{0.758 cm}$\displaystyle 0<16\,\ep_o<
\varlimsup_{j\to\infty}\;\sup \Big\{\, \frac 1{b_j-a_j+1}\,
\sum\limits_{k=a_j}^{b_j}\, |\langle x^*, x_{\substack{{}\\ k}}
\rangle |\, ;\, x^*\in X^* ,\,\|x^*\|\leq 1\,\Big\}\, .$
\smallskip

\noindent Then
\smallskip

\noindent\hspace{0.281 cm}$\displaystyle
\mathcal J := \bigg\{\, j\geq 1\, ;\, \sup \Big\{\,
\frac 1{b_j-a_j+1}\,\sum\limits_{k=a_j}^{b_j}\, |\langle x^*,
x_{\substack{{}\\ k}}\rangle |\, ;\, x^*\in X^* ,\,\|x^*\|
\leq 1\,\Big\} > 16\,\ep_o\,\bigg\}$
\smallskip

\noindent is infinite. Using (in the complex case)
$\langle x^*,x_{\substack{{}\\ k}}\rangle =\Re\,\langle\, x^*,
x_{\substack{{}\\ k}}\rangle -i\,\Re\,\langle\, i\, x^*,
x_{\substack{{}\\ k}}\rangle\, ,$
it follows that also
\smallskip

\noindent\hspace{0.106 cm}$\displaystyle
\mathcal J_{\substack{{}\\ \Re}} := \bigg\{\, j\geq 1\, ;\,
\sup \Big\{\,\frac 1{b_j-a_j+1}\,\sum\limits_{k=a_j}^{b_j}\,
|\,\Re\,\langle x^*,x_{\substack{{}\\ k}}\rangle |\, ;\, x^*
\in X^* ,\,\|x^*\|\leq 1\,\Big\} > 8\,\ep_o\,\bigg\}$
\smallskip

\noindent is infinite. Now, since $\,\Re\,\langle x^*,
x_{\substack{{}\\ k}}\rangle =\Re^+\langle\, x^*,
x_{\substack{{}\\ k}}\rangle -\Re^+\langle\, -\, x^*,
x_{\substack{{}\\ k}}\rangle\,$, we obtain that
\smallskip

\noindent\hspace{0.139 cm}$\displaystyle
\mathcal J_{\substack{{}\\ +}} := \bigg\{\, j\geq 1\, ;\,
\sup \Big\{\,\frac 1{b_j-a_j+1}\,\sum\limits_{k=a_j}^{b_j}\,
\Re^+\langle x^*,x_{\substack{{}\\ k}}\rangle\, ;\, x^*
\in X^* ,\,\|x^*\|\leq 1\,\Big\} > 4\,\ep_o\,\bigg\}$
\smallskip

\noindent is infinite.

Let $j\in\mathcal J_{\substack{{}\\ +}}$ be arbitrary. Then
there exists $y_j^*\in X^*$ with $\| y_j^*\|\leq 1$ such that
\smallskip

\noindent\hspace{3.566 cm}$\displaystyle \frac 1{b_j-a_j+1}\,
\sum\limits_{k=a_j}^{b_j}\,\Re^+\langle\, y_j^*,
x_{\substack{{}\\ k}}\rangle > 4\,\ep_o\, .$
\smallskip

\noindent Denoting $\,\mathcal B_j:=\{ a_j\leq k\leq b_j\, ;\,
\Re^+\langle\, y_j^*,x_{\substack{{}\\ k}}\rangle >2\,\ep_o\}\,$,
we have
\medskip

\noindent\hspace{1.863 cm}$\displaystyle 4\,\ep_o <\frac 1{b_j-a_j+1}
\,\Big(\sum\limits_{k\in\mathcal B_j}\Re^+\langle\, y_j^*,
x_{\substack{{}\\ k}}\rangle +\!\sum_{\substack{a_j\leq k\leq b_j \\
k\notin\mathcal B_j}}\!\Re^+\langle\, y_j^*,
x_{\substack{{}\\ k}}\rangle\Big)$

\noindent\hspace{2.514 cm}$\displaystyle \leq \frac 1{b_j-a_j+1}\,
\mathrm{card}\, (\mathcal B_j ) +2\,\ep_o\, ,$
\medskip

\noindent hence $\,\mathrm{card}\, (\mathcal B_j )\geq 2\,
(b_j-a_j+1)\,\ep_o\,$.

Denoting now by $j_{\substack{{}\\ 1}}$ the least element of
$\mathcal J_{\substack{{}\\ +}}\,$, we can construct recursively
a sequnce $j_{\substack{{}\\ 1}} < j_{\substack{{}\\ 2}} <\,
\ldots\,$ in $\mathcal J_{\substack{{}\\ +}}$ such that, for every
$n\geq 2\,$,
\smallskip

\noindent\hspace{4.788 cm}$b_{j_n}-b_{j_{n-1}} > n\,$ and
\smallskip

\noindent\hspace{0.755 cm}the cardinality of
$\,{\mathcal B_{\substack{{}\\ j_n}}}' :=
\{\, k\in\mathcal B_{\substack{{}\\ j_n}}\, ;\, k>b_{j_{n-1}}+n
\,\}\,$ is $\, > (b_{j_n}-a_{j_n})\,\ep_o\,$.
\smallskip

\noindent Indeed, if we choose $j_n$ in the infinite set
$\mathcal J_{\substack{{}\\ +}}$ such that $j_n>j_{n-1}$ and
\smallskip

\noindent\hspace{2.663 cm}$\displaystyle b_{j_n}-a_{j_n}+1\geq
j_n +1 >\frac{b_{j_{n-1}}+n}{\ep_o}\geq b_{j_{n-1}}+n\,$,
\smallskip

\noindent then $b_{j_n}-b_{j_{n-1}} > n$ and
\smallskip

\noindent\hspace{0.728 cm}$\mathrm{card}\, (
{\mathcal B_{\substack{{}\\ j_n}}}') \geq \mathrm{card}\, (
\mathcal B_{\substack{{}\\ j_n}}) - (b_{j_{n-1}}+n) \geq
2\, (b_{j_n}-a_{j_n}+1)\,\ep_o - (b_{j_{n-1}}+n)$
\smallskip

\noindent\hspace{2.46 cm}$ >(b_{j_n}-a_{j_n}+1)\,\ep_o\, .$
\medskip

Putting
\smallskip

\noindent\hspace{5.23 cm}$\displaystyle \mathcal B :=
\bigcup_{n\geq 2}\, {\mathcal B_{\substack{{}\\ j_n}}}'\, ,$
\smallskip

\noindent we have for every $n\geq 2$
\smallskip

\noindent\hspace{1.143 cm}$\mathcal B\cap (\, b_{j_{n-1}}\, ,\,
b_{j_{n-1}}+n\, ]=\emptyset\, ,\qquad
\mathcal B\cap (\, b_{j_{n-1}}+n\, ,\, b_{j_n}\, ] =
{\mathcal B_{\substack{{}\\ k_j}}}' \subset
\mathcal B_{\substack{{}\\ k_j}}\, ,$
\smallskip

\noindent and so

\noindent\hspace{1.26 cm}$\Re\,\langle\, y_{j_n}^{\, *},
x_{\substack{{}\\ k}}\rangle =\Re^+\langle\, y_{j_n}^{\, *},
x_{\substack{{}\\ k}}\rangle > 2\,\ep_o\,$ for all $\, k\in
\mathcal B\cap (\, b_{j_{n-1}}+n\, ,\, b_{j_n}\, ]\, .$
\medskip

\noindent On the other hand,
\medskip

\noindent\hspace{1.8 cm}$\displaystyle BD^*(\mathcal B) \geq
\varlimsup_{n\to\infty}\,\frac 1{b_{j_n}-a_{j_n}+1}\;
\mathrm{card}\,\big( \mathcal B \cap [a_{j_n},b_{j_n}]\big)$
\smallskip

\noindent\hspace{3.158 cm}$\displaystyle \geq
\varlimsup_{n\to\infty}\,\frac 1{b_{j_n}-a_{j_n}+1}\;
\mathrm{card}\,\big( \underbrace{\mathcal B_{\substack{{}\\ j_n}}'
\cap [a_{j_n},b_{j_n}]}_{=\,\mathcal B_{\substack{{}\\ j_n}}'}
\big)\geq\ep_o\, .$

\noindent Therefore, setting $x_n^{\;\! *}:=y_{j_n}^{\, *}\,$,
the proof is complete.

\end{proof}

For the proof of the next theorem we adapt the proof of
Proposition \ref{refl_wm=unif.wm}, in which instead of
Lemmas \ref{HB}, \ref{Jones} and \ref{weak-structure} we use
Lemma \ref{BHB}, Theorem \ref{structure} and Lemma
\ref{zeroinclosure} :

\begin{theorem}[Weak mixing for convex shift-bounded
sequences]\label{wm=unif.wm}
For a convex shift-bounded sequence
$(x_{\substack{{}\\ k}})_{\substack{{}\\ k\geq 1}}$ in a
Banach space $X\,$, the following conditions are equivalent $:$
\begin{itemize}
\item[(i)] $(x_{\substack{{}\\ k}})_{\substack{{}\\ k\geq 1}}$ is
weakly mixing to zero, that is

\begin{center} $\displaystyle \lim_{n\to\infty}\; \frac 1n\,
\sum_{k=1}^n\, |\langle x^*,x_{\substack{{}\\ k}}\rangle | = 0\;\,
\text{for all}\;\, x^*\in X^*\,$.
\end{center}
\item[(j)] $(x_{\substack{{}\\ k}})_{\substack{{}\\ k\geq 1}}$ is
uniformly weakly mixing to zero, that is

\begin{center} $\displaystyle \lim_{n\to\infty}\;
\sup\, \Big\{\, \frac 1n\, \sum\limits_{k=1}^n\, |\langle x^*,
x_{\substack{{}\\ k}}\rangle |\, ;\, x^*\in X^*\, ,\,\|x^*\|
\leq 1\,\Big\}= 0\,$.
\end{center}
\item[(jw)] {\rm (\ref{B.unif.w.mix.})} holds, that is
\smallskip

\noindent\hspace{0.741 cm}$\displaystyle \varlimsup_{\substack{a,b
\in\mathbb N^*\\ b-a\to\infty}}\;\sup \Big\{\, \frac 1{b-a+1}\,
\sum\limits_{k=a}^b\, |\langle x^*, x_{\substack{{}\\ k}}\rangle
|\, ;\, x^*\in X^*\, ,\,\|x^*\|\leq 1\,\Big\}= 0\,$.
\end{itemize}
\end{theorem}

\begin{proof}
The implications (jw)$\,\Rightarrow\,$(j)$\,\Rightarrow\,$(i)
are trivial. For (i)$\,\Rightarrow\,$(jw) we shall show that
(i) and the negation of (jw) lead to a contradiction.

Let $c>0$ be a constant such that (\ref{pow.bdd}) holds for
any choice of $p\in\mathbb{N}^*$ and $\lambda_1\, ,\,\ldots\,
,\,\lambda_p\geq 0\,$. Since (\ref{B.unif.w.mix.}) does not hold,
there exist $a_{\substack{{}\\ 1}}\, ,\, b_{\substack{{}\\ 1}}\, ,\,
a_{\substack{{}\\ 2}}\, ,\, b_{\substack{{}\\ 2}}\, ,\,\ldots\,
\in\mathbb N^*\,$ with $\, b_j-a_j\geq j\, ,\, j\geq 1\,$, such
that
\smallskip

\noindent\hspace{1.326 cm}$\displaystyle \varlimsup_{j\to\infty}
\;\sup \Big\{\, \frac 1{b_j-a_j+1}\, \sum\limits_{k=a_j}^{b_j}
\, |\langle x^*, x_{\substack{{}\\ k}}\rangle |\, ;\, x^*\in
X^*\, ,\,\|x^*\|\leq 1\,\Big\} >0\, .$
\smallskip

By Lemma \ref{BHB} there exist
\smallskip

\noindent\hspace{5.442 cm}$0<\ep_o\leq 1\,$,

\noindent\hspace{4.173 cm}$\mathcal B\subset\mathbb{N}^*$ with
$\, BD^*(\mathcal B)\geq\ep_o\,$,

\noindent\hspace{3.135 cm}$j_{\substack{{}\\ 1}} <
j_{\substack{{}\\ 2}} <\, \ldots\,$ in $\,\mathbb{N}^*\!$ with
$\, b_{j_n}-b_{j_{n-1}} > n\,$,

\noindent\hspace{3.741 cm}$x_1^*\, ,\, x_2^*\, ,\, \ldots\,\in
X^*$ with $\,\| x_n^*\|\leq 1\,$,
\smallskip

\noindent for which
\smallskip

\noindent\hspace{3.943 cm}$\displaystyle \mathcal B \cap
\bigcup_{n\geq 2} (\, b_{j_{n-1}}\, ,\, b_{j_{n-1}}+n\, ]
=\emptyset\, ,$

\noindent\hspace{2.037 cm}$\Re\,\langle\, x_n^*,
x_{\substack{{}\\ k}}\rangle > 2\,\ep_o\, ,\qquad k\in
\mathcal B\cap (\, b_{j_{n-1}}+n\, ,\, b_{j_n}\, ]\; ,\, n
\geq 2\, .$
\medskip

\noindent Further, by Theorem \ref{structure}, there exist
\medskip

\noindent\hspace{3.627 cm}$\mathcal A\subset\mathbb N^*\,$
having density $\, D(\mathcal A) >0\, ,$

\noindent\hspace{2.365 cm}$1\leq m_{\substack{{}\\ 1}}<
m_{\substack{{}\\ 2}}<\,\ldots\;$ and $\, 1\leq
n_{\substack{{}\\ 1}}<n_{\substack{{}\\ 2}}<\,\ldots\;$ in $\,
\mathbb N^* ,$
\smallskip

\noindent such that
\smallskip

\noindent\hspace{2.332 cm}$\mathcal A\cap [1,m_j] =\big\{ k\in
[1,m_j]\, ;\, k+n_j\in\mathcal B\big\}\; ,\qquad j\geq 1\,$.
\smallskip

\noindent Finally, (i) and Lemma \ref{zeroinclosure} entail
that there are $p\in\mathbb{N}^*$, $k_{\substack{{}\\ 1}}<\,
\ldots\, <k_p$ in $\mathcal A\,$ and $\lambda_{\substack{{}\\ 1}}\,
,\,\ldots\, ,\,\lambda_p\geq 0\, ,\,\lambda_{\substack{{}\\ 1}}
+ \ldots +\lambda_p=1\,$, such that
\smallskip

\noindent\hspace{4.808 cm}$\displaystyle \Big\|\sum_{j=1}^p\,
\lambda_j\, x_{\substack{{}\\ k_j}}\Big\|\leq\frac{\ep_o}c\, .$
\medskip

\noindent By (\ref{pow.bdd}) it follows that
%
\begin{equation}\label{transl2}
\Big\|\sum_{j=1}^p\,\lambda_j\, x_{\substack{{}\\ {k_j +n}}}
\Big\|\leq c\,\Big\|\sum_{j=1}^p\,\lambda_j\,
x_{\substack{{}\\ k_j}}\Big\|\leq\ep_o\; ,\qquad n\geq 1\, .
\end{equation}

Now let $q\in\mathbb N^*$ be such that $k_{\substack{{}\\ 1}}
\, ,\,\ldots\, ,\, k_p\leq m_q\,$. Then
\smallskip

\noindent\hspace{3.415 cm}$k_{\substack{{}\\ 1}}+n_j\, ,\,
\ldots\, ,\, k_p+n_j\in\mathcal B\, ,\qquad j\geq q\, .$
\smallskip

\noindent Choose $j_*\geq q$ with $k_{\substack{{}\\ 1}}
+n_{j_*}\geq b_{j_{m_q}}$ and define $n\in\mathbb N^*$ by
$b_{j_{n-1}}<k_{\substack{{}\\ 1}}+n_{j_*}\leq b_{j_n}\,$.
Since $b_{j_{m_q}}\leq k_{\substack{{}\\ 1}}+n_{j_*}\leq
b_{j_n}$ and the sequence $(b_{j_{n'}})_{\substack{{}\\ n'\geq 1}}$
is increasing, we have $m_q\leq n\,$. We claim that
\begin{equation}\label{lacunary2}
b_{j_{n-1}}+n<k_{\substack{{}\\ 1}}+n_{j_*}\leq k_p+n_{j_*}\leq
b_{j_n}\, .
\end{equation}

Indeed, $b_{j_{n-1}}<k_{\substack{{}\\ 1}}+n_{j_*}\,$,
$k_{\substack{{}\\ 1}}+n_{j_*}\in\mathcal B\,$ and $\,\mathcal B
\cap (b_{j_{n-1}} , b_{j_{n-1}}+n]=\emptyset$ imply
$b_{j_{n-1}}+n <k_{\substack{{}\\ 1}}+n_{j_*}\,$. Further,
$k_p+n_{j_*} =k_{\substack{{}\\ 1}}+n_{j_*} +(k_p-
k_{\substack{{}\\ 1}})\leq b_{j_n}+m_q <b_{j_n}+n+1\,$,
$k_p+n_{j_*}\in\mathcal B\,$ and $\,\mathcal B\cap
(b_{j_n},b_{j_n}+n+1]=\emptyset\,$ yield $\, k_p+n_{j_*}\leq
b_{j_n}\,$.

By (\ref{lacunary2}) we have $\, k_{\substack{{}\\ 1}}+n_{j_*}
\, ,\,\ldots\, ,\, k_p+n_{j_*}\in\mathcal B\cap (b_{j_{n-1}}+n\,
,b_{j_n}]\,$, so
\smallskip

\noindent\hspace{3.375 cm}$\Re\,\langle\, x_n^*,
x_{\substack{{}\\ k_j+n_{j_*}}}\rangle > 2\,\ep_o\, ,\qquad
1\leq j\leq p\, .$
\smallskip

\noindent Since $\| x_n^*\|\leq 1\,$, it follows that
\smallskip

\noindent\hspace{0.653 cm}$\displaystyle \Big\|\sum_{j=1}^p\,
\lambda_j\, x_{\substack{{}\\ {k_j +n}}}\Big\|\geq \Re\,\left<
x_n^*,\,\sum_{j=1}^p\,\lambda_j\, x_{\substack{{}\\ {k_j +n}}}
\right> =\sum_{j=1}^p\,\lambda_j\,\Re\,\langle\, x_n^*,
x_{\substack{{}\\ k_j+n_{j_*}}}\rangle > 2\,\ep_o\, ,$
\smallskip

\noindent in contradiction with (\ref{transl2}).

\end{proof}

\bigskip
\section{Weak mixing to zero for Cesaro shift-bounded sequences}
\label{Cesaro}

If $X$ is a uniformly convex Banach space, then Theorem
\ref{wm=unif.wm} holds under a milder assumption on
$(x_{\substack{{}\\ k}})_{\substack{{}\\ k\geq 1}}$ than convex
shift-boundedness.

For we shall use that in uniformly convex Banach spaces the
classical Mazur theorem about the equality of the weak and norm
closure of a convex subset holds in the following sharper form :

\begin{theorem}[Mazur type theorem in uniformly convex Banach
spaces]\label{Mazur}
Let $S$ be a bounded subset of a uniformly convex Banach space
$X\,$, and $x$ an element of the weak closure of $S\,$. Then
there exists a sequence
$(x_{\substack{{}\\ k}})_{\substack{{}\\ k\geq 1}}\subset S\,$,
such that
\smallskip

\noindent\hspace{4.18 cm}$\displaystyle \lim_{n\to\infty}\;
\Big\|\, x-\frac 1n\,\sum_{k=1}^n\, x_{\substack{{}\\ k}}\,
\Big\| =0\;$.
\end{theorem}

\begin{proof}
Uniformly convex Banach spaces are reflexive (see e.g. \cite{D},
page 131), so $S$ is weakly relatively compact.
Consequently, since normed linear spaces are angelic in their
weak topology (see e.g. \cite{Fl}, 3.10.(1)), there exists a
sequence $(y_{\substack{{}\\ j}})_{\substack{{}\\ j\geq 1}}$
in $S\,$, which is weakly convergent to $x\,$. Now, according
to the validity of the ``Banach-Saks Theorem'' \cite{BaS} in
uniformly convex Banach spaces, due to S. Kakutani \cite{Ka}
(see also \cite{D}, Chapter VIII, Theorem 1), there exists a
subsequence $\big( y_{j_{\substack{{}\\ k}}}
\big)_{\substack{{}\\ k\geq 1}}$ such that
\smallskip

\noindent\hspace{4.13 cm}$\displaystyle \lim_{n\to\infty}\;
\Big\|\, x-\frac 1n\,\sum_{k=1}^n\, y_{j_{\substack{{}\\ k}}}\,
\Big\| =0\;$.

\end{proof}

Let us call a sequence
$(x_{\substack{{}\\ k}})_{\substack{{}\\ k\geq 1}}$ in a Banach
space $X$ \emph{Cesaro shift-bounded} if there exists a constant
$c>0$ such that (\ref{pow.bdd}) holds for any choice of
$p\in\mathbb{N}^*$ and $\,\lambda_1\, ,\,\ldots\, ,\,\lambda_p
\in\{ 0,1\}\,$, that is
\begin{center}
$\displaystyle \Big\|\sum_{j=1}^p
\, x_{n_j+n}\,\Big\|\leq c\,\Big\|
\sum_{j=1}^p\, x_{n_j}\,\Big\|\, ,\qquad n\geq 1$
\end{center}
for any $p\in\mathbb N^*$ and $\, n_1\, ,\,\ldots\, ,\, n_p\in
\mathbb N^*$ with $\, n_1<\,\ldots\, <n_p\,$.

Clearly, every convex shift-bounded sequence in
$X$ is Cesaro shift-bounded, but the converse does not hold,
even in Hilbert spaces :

\begin{example}\label{ex4}
{\it Let $H$ be an infinite-dimensional Hilbert space and
choose, for every $k\in\mathbb N$, three vectors
$u_{\substack{{}\\ k}}\, ,\, v_{\substack{{}\\ k}}\, ,\,
w_{\substack{{}\\ k}}\in H$ such that
\medskip

\noindent\hspace{0.133 cm}$\displaystyle \| u_{\substack{{}\\ k}}\|
=\| v_{\substack{{}\\ k}} \| =\| w_{\substack{{}\\ k}}\| =1\, ,\,
0<\| u_{\substack{{}\\ k}}-v_{\substack{{}\\ k}}\| <\frac 1{k+3}
\, ,\, w_{\substack{{}\\ k}}\perp \{u_{\substack{{}\\ k}}\, ,\,
v_{\substack{{}\\ k}}\}\,$ for all $\,\;\! k\in\mathbb N$,

\noindent\hspace{2.773 cm}$\{u_{\substack{{}\\ k}}\, ,\,
v_{\substack{{}\\ k}}\, ,\, w_{\substack{{}\\ k}}\}\perp
\{u_{\substack{{}\\ l}}\, ,\, v_{\substack{{}\\ l}}\, ,\,
w_{\substack{{}\\ l}}\}\,$ whenever $\,\;\! k\neq l\,$.
\medskip

\noindent Let us define the sequence $(x_n)_{n\geq 1}$ by
\medskip

\noindent\hspace{1.746 cm}$x_{\substack{{}\\ 3k+1}}:=2\,
u_{\substack{{}\\ k}}\, ,\, x_{\substack{{}\\ 3k+2}}:=
-v_{\substack{{}\\ k}}\, ,\, x_{\substack{{}\\ 3k+3}}:=
w_{\substack{{}\\ k}}\,$ for even $\, k\in\mathbb N\,$,

\noindent\hspace{1.746 cm}$x_{\substack{{}\\ 3k+1}}:=2\,
u_{\substack{{}\\ k}}\, ,\, x_{\substack{{}\\ 3k+2}}:=
w_{\substack{{}\\ k}}\, ,\, x_{\substack{{}\\ 3k+3}}:=
-v_{\substack{{}\\ k}}\,$ for odd $\, k\in\mathbb N\,$.
\medskip

\noindent Then $(x_n)_{n\geq 1}$ is Cesaro shift-bounded,
but not convex shift-bounded.}
\end{example}

\begin{proof}
For every $k\in\mathbb N$ we denote by
$\mathcal V_{\substack{{}\\ k}}$ the set of the vectors
\medskip

\noindent\hspace{4.744 cm}$x_{\substack{{}\\ 3k+1}}\, ,\,
x_{\substack{{}\\ 3k+2}}\, ,\, x_{\substack{{}\\ 3k+3}}$

\noindent\hspace{2.8 cm}$x_{\substack{{}\\ 3k+1}}+
x_{\substack{{}\\ 3k+2}}\, ,\,x_{\substack{{}\\ 3k+1}}+
x_{\substack{{}\\ 3k+3}}\, ,\,x_{\substack{{}\\ 3k+2}}+
x_{\substack{{}\\ 3k+3}}$

\noindent\hspace{4.513 cm}$x_{\substack{{}\\ 3k+1}}+
x_{\substack{{}\\ 3k+2}}+x_{\substack{{}\\ 3k+3}}\,$.
\medskip

\noindent It is easy to verify that $\displaystyle \frac 23\leq
\| x\|\leq\sqrt{5}\,$ for all $\, x\in\mathcal V_k\,$.

Now let $p\in\mathbb N^*$, $\, n_1\, ,\,\ldots\, ,\, n_p\in
\mathbb N^*$ with $\, n_1<\,\ldots\, <n_p\,$, and $n\in
\mathbb N^*$ be arbitrary. Let $q$
denote the number of all $\mathcal V_{\substack{{}\\ k}}$
which contain some $x_{n_j}\,$. Then
\begin{equation}\label{greater}
\Big\|\sum_{j=1}^p\, x_{n_j}\,\Big\|\geq\sqrt{q\,
\Big(\frac 23\Big)^2} =\sqrt{\frac 49\,q}\, .
\end{equation}
On the other hand, since the number of all
$\mathcal V_{\substack{{}\\ k}}$ which contain some
$x_{n_j +n}$ is $\,\leq 2\, q\,$, we have
\begin{equation}\label{less}
\Big\|\sum_{j=1}^p\, x_{n_j+n}\,\Big\|\leq\sqrt{2\, q\,
\big(\sqrt{5}\big)^2} =\sqrt{10\, q}\, .
\end{equation}
Now (\ref{less}) and (\ref{greater}) entail that
\smallskip

\noindent\hspace{3.821 cm}$\displaystyle \Big\|\sum_{j=1}^p
\, x_{n_j+n}\,\Big\|\leq\sqrt{\frac{45}2}\,\Big\|
\sum_{j=1}^p\, x_{n_j}\,\Big\|$
\smallskip

\noindent and we conclude that the sequence $(x_n)_{n\geq 1}$
is Cesaro shift-bounded.

To show that $(x_n)_{n\geq 1}$ is not convex shift-bounded,
let us assume the contrary, that is the existence of some
constant $c>0$ such that
\medskip

\noindent\hspace{2.82 cm}$\displaystyle \Big\|\sum_{n=1}^p\,
\lambda_n\, x_{n+m}\,\Big\|\leq c\,\Big\|\sum_{n=1}^p\,
\lambda_n\, x_n\,\Big\|\, ,\qquad m\geq 1$
\smallskip

\noindent for any $p\in\mathbb N^*$ and $\lambda_1\, ,\,
\ldots\, ,\,\lambda_p\geq 0\,$. Let $k$ be an arbitrary
even number in $\mathbb N^*$ and set
$p:=3k+2\,$, $\lambda_n:=0$ for $1\leq n\leq 3k\,$,
$\lambda_{3k+1}:=1$ and $\lambda_{3k+2}:=2\,$. Then
\begin{center}
$\displaystyle \Big\|\sum_{n=1}^p\,
\lambda_n\, x_n\,\Big\| =\big\|\, 2\, u_{\substack{{}\\ k}}
-2\, v_{\substack{{}\\ k}}\big\| <\frac 2{k+3}\,$,
\end{center}
while
\begin{center}
$\displaystyle \Big\|\sum_{n=1}^p\,
\lambda_n\, x_{n+3}\,\Big\| =\big\|\, 2\,
u_{\substack{{}\\ k+1}} +2\, w_{\substack{{}\\ k+1}}\big\|
=2\,\sqrt{2}\,$.
\end{center}
It follows that $\,\displaystyle 2\,\sqrt{2}\leq
\frac{2\, c}{k+3}\,$, what is not possible for any even
$k\in\mathbb N^*$.

\end{proof}

Using Theorem \ref{Mazur}, we can adapt the proof of Theorem
\ref{wm=unif.wm} to the case of Cesaro shift-bounded sequences
in uniformly convex Banach spaces :

\begin{theorem}[Weak mixing for Cesaro shift-bounded
sequences]\label{Cesaro-bdd}
For a Cesaro shift-bounded sequence
$(x_{\substack{{}\\ k}})_{\substack{{}\\ k\geq 1}}$ in a
uniformly convex Banach space $X\,$, the following conditions
are equivalent $:$
\begin{itemize}
\item[(i)] $(x_{\substack{{}\\ k}})_{\substack{{}\\ k\geq 1}}$ is
weakly mixing to zero,
\item[(j)] $(x_{\substack{{}\\ k}})_{\substack{{}\\ k\geq 1}}$ is
uniformly weakly mixing to zero,
\item[(jw)] {\rm (\ref{B.unif.w.mix.})} holds.
\end{itemize}

\end{theorem}

\begin{proof}
Since the implications (jw)$\,\Rightarrow\,$(j)$\,\Rightarrow\,$(i)
are trivial, to complete the proof we need only to prove that
(i)$\,\Rightarrow\,$(jw). 

Let $c>0$ be a constant such that (\ref{pow.bdd}) holds for
any choice of $p\in\mathbb{N}^*$ and $\lambda_1\, ,\,\ldots\,
,\,\lambda_p\in\{ 0,1\}\,$. Since (\ref{B.unif.w.mix.}) does not hold,
there exist $a_{\substack{{}\\ 1}}\, ,\, b_{\substack{{}\\ 1}}\, ,\,
a_{\substack{{}\\ 2}}\, ,\, b_{\substack{{}\\ 2}}\, ,\,\ldots\,
\in\mathbb N^*\,$ with $\, b_j-a_j\geq j\, ,\, j\geq 1\,$, such
that
\smallskip

\noindent\hspace{1.326 cm}$\displaystyle \varlimsup_{j\to\infty}
\;\sup \Big\{\, \frac 1{b_j-a_j+1}\, \sum\limits_{k=a_j}^{b_j}
\, |\langle x^*, x_{\substack{{}\\ k}}\rangle |\, ;\, x^*\in
X^*\, ,\,\|x^*\|\leq 1\,\Big\} >0\, .$
\smallskip

By Lemma \ref{BHB} there exist
\smallskip

\noindent\hspace{5.442 cm}$0<\ep_o\leq 1\,$,

\noindent\hspace{4.173 cm}$\mathcal B\subset\mathbb{N}^*$ with
$\, BD^*(\mathcal B)\geq\ep_o\,$,

\noindent\hspace{3.135 cm}$j_{\substack{{}\\ 1}} <
j_{\substack{{}\\ 2}} <\, \ldots\,$ in $\,\mathbb{N}^*\!$ with
$\, b_{j_n}-b_{j_{n-1}} > n\,$,

\noindent\hspace{3.741 cm}$x_1^*\, ,\, x_2^*\, ,\, \ldots\,\in
X^*$ with $\,\| x_n^*\|\leq 1\,$,
\smallskip

\noindent for which
\smallskip

\noindent\hspace{3.943 cm}$\displaystyle \mathcal B \cap
\bigcup_{n\geq 2} (\, b_{j_{n-1}}\, ,\, b_{j_{n-1}}+n\, ]
=\emptyset\, ,$

\noindent\hspace{2.037 cm}$\Re\,\langle\, x_n^*,
x_{\substack{{}\\ k}}\rangle > 2\,\ep_o\, ,\qquad k\in
\mathcal B\cap (\, b_{j_{n-1}}+n\, ,\, b_{j_n}\, ]\; ,\, n
\geq 2\, .$
\medskip

On the other hand, since $X$ is reflexive and any bounded set
in a reflexive Banach space is weakly relatively compact, by
Lemma \ref{Jones} there exists $\mathcal A_o\subset\mathbb N^*$
with $D_*(\mathcal A_o) =1$ such that $\displaystyle
\lim_{\mathcal{A}_o\ni k\to\infty} x_k =0\,$ in the weak
topology of $X\,$.

Finally, by Corollary \ref{suff.-structure} there exists
$\mathcal I\subset \mathcal A_o$ with $D^*(\mathcal I)>0\,$,
such that
\medskip

\noindent\hspace{2.214 cm}$\{ n\in\mathbb N\, ;\,\mathcal F +n
\subset\mathcal B\}\,$ is infinite for any finite $\,
\mathcal F\subset\mathcal I\, .$
\medskip

\noindent Then $\displaystyle \lim_{\mathcal{I}\ni k\to\infty}
x_k =0\,$ with respect to the weak topology of $X$ and by
Theorem \ref{Mazur} there are $p\in\mathbb{N}^*$ and
$k_{\substack{{}\\ 1}}<\,\ldots\, <k_p$ in $\mathcal I\,$
such that
\medskip

\noindent\hspace{4.81 cm}$\displaystyle \Big\|\,\frac 1p\,
\sum_{j=1}^p\, x_{\substack{{}\\ k_j}}\,\Big\|\leq
\frac{\varepsilon_o}c\,$.
\smallskip

\noindent By (\ref{pow.bdd}) it follows that
\begin{equation}\label{transl3}
\Big\|\,\frac 1p\,\sum_{j=1}^p\, x_{\substack{{}\\ k_j +n}}\,
\Big\|\leq c\,\Big\|\,\frac 1p\,\sum_{j=1}^p\,
x_{\substack{{}\\ k_j}}\,\Big\|\leq\ep_o\; ,\qquad n\geq 1\, .
\end{equation}

Now set $m:=\max\big( k_p-k_{\substack{{}\\ 1}},2\big)\,$.
Since the set $\big\{ k\in\mathbb N^*\, ;\,\{ k_{\substack{{}\\ 1}}
\, ,\,\ldots\, ,\, k_p\} +n\subset\mathcal B\big\}$ is infinite,
it contains some $n_o$ such that $k_{\substack{{}\\ 1}} +n_o
\geq b_{j_m}\,$. We define $n_{\substack{{}\\ 1}}\in\mathbb N^*$
by $b_{j_{n_1 -1}}<k_{\substack{{}\\ 1}}+n_o\leq b_{j_{n_1}}\,$.
Since $b_{j_m}\leq k_{\substack{{}\\ 1}}+n_o\leq b_{j_{n_1}}$
and the sequence $(b_{j_n})_{\substack{{}\\ n\geq 1}}$ is
increasing, we have $m\leq n_{\substack{{}\\ 1}}\,$. We claim that
\begin{equation}\label{lacunary3}
b_{j_{n_1 -1}}+n_{\substack{{}\\ 1}}\leq k_{\substack{{}\\ 1}}+n_o
<k_p+n_o\leq b_{j_{n_1}}\, .
\end{equation}

Indeed, $b_{j_{n_1 -1}}<k_{\substack{{}\\ 1}}+n_o\,$,
$k_{\substack{{}\\ 1}}+n_o\in\mathcal B\,$ and $\,\mathcal B
\cap (b_{j_{n_1 -1}} , b_{j_{n_1 -1}}+n_{\substack{{}\\ 1}}]=
\emptyset$ imply $b_{j_{n_1 -1}}+n_{\substack{{}\\ 1}}
<k_{\substack{{}\\ 1}}+n_o\,$. Further, $k_p+n_o =
k_{\substack{{}\\ 1}}+n_o +(k_p-k_{\substack{{}\\ 1}})\leq
b_{j_{n_1}}+m <b_{j_{n_1}}+n_{\substack{{}\\ 1}}+1\,$,
$k_p+n_o\in\mathcal B\,$ and $\,\mathcal B\cap
(b_{j_{n_1}},b_{j_{n_1}}+n_{\substack{{}\\ 1}}+1]=\emptyset\,$
yield $\, k_p+n_o\leq b_{j_{n_1}}\,$.

By (\ref{lacunary3}) we have $\, k_{\substack{{}\\ 1}}+n_o\, ,\,
\ldots\, ,\, k_p+n_o\in\mathcal B\cap (b_{j_{n_1 -1}}+
n_{\substack{{}\\ 1}}\, ,b_{j_{n_1}}]\,$, so
\smallskip

\noindent\hspace{3.372 cm}$\Re\,\langle\, x_{n_1}^*,
x_{\substack{{}\\ k_j+n_o}}\rangle > 2\,\ep_o\, ,\qquad
1\leq j\leq p\, .$
\smallskip

\noindent Since $\| x_{n_1}^*\|\leq 1\,$, it follows that
\smallskip

\noindent\hspace{0.43 cm}$\displaystyle \Big\|\,\frac 1p\,
\sum_{j=1}^p\, x_{\substack{{}\\ k_j +n_o}}\,\Big\|\geq \Re\,
\left< x_{n_1}^*,\,\frac 1p\,\sum_{j=1}^p\,
x_{\substack{{}\\ k_j +n_o}}\,\right> =\frac 1p\,\sum_{j=1}^p\,
\Re\,\langle\, x_{n_1}^*, x_{\substack{{}\\ k_j+n_o}}\rangle >
2\,\ep_o\, ,$
\smallskip

\noindent in contradiction with (\ref{transl3}).

\bigskip

\end{proof}

\bigskip
\section{Appendix: Ergodic theorem for convex shift-bounded
sequences}

We can define also \emph{ergodicity} of a bounded sequence
$(x_{\substack{{}\\ k}})_{\substack{{}\\ k\geq 1}}$ in a Banach
space $X$ by requiring that
\begin{center}
$\displaystyle \lim_{n\to\infty}\; \Big\|\,\frac 1n\,\sum_{k=1}^n
x_{\substack{{}\\ k}}\Big\| =0$
\end{center}
(cf. \cite{BB}, Section 3). Clearly, if
$(x_{\substack{{}\\ k}})_{\substack{{}\\ k\geq 1}}$ is uniformly
weak mixing, then it is ergodic. In this section we complete our
knowledge about convex shift-bounded sequences by proving a
mean ergodic theorem for them (Corollary \ref{pow.bdd.m.erg.}).

Let $l^\infty (X)$ denote the vector space of all bounded
sequnces $(x_{\substack{{}\\ j}})_{\substack{{}\\ j\geq 1}}$ in a
Banach space $X\,$, endowed with the uniform norm
$\big\| (x_{\substack{{}\\ j}})_{\substack{{}\\ j}}\big\|_\infty =
\sup_{j} \| x_{\substack{{}\\ j}}\|\,$, and let
$\sigma_{\substack{{}\\ \!\leftarrow}}$ be the backward
shift on $l^\infty (X)\,$, defined by
\smallskip

\noindent\hspace{4.27 cm}$\sigma_{\substack{{}\\ \!\leftarrow}}
\big( (x_{\substack{{}\\ j}})_{\substack{{}\\ j\geq 1}}\big) =
(x_{\substack{{}\\ j+1}})_{\substack{{}\\ j\geq 1}}\,$.

\begin{theorem}[Mean Ergodic Theorem for sequences]\label{m.erg.}
For a bounded sequence
$(x_{\substack{{}\\ k}})_{\substack{{}\\ k\geq 1}}$ in a
Banach space $X\,$, the following conditions are equivalent $:$
\smallskip

\noindent\hspace{0.62 cm}{\rm (e)}\hspace{3.1 cm}
$\displaystyle \lim_{\substack{0\leq m<n\\ n-m\to\infty}}\; \Big\|
\,\frac 1{n-m}\,\sum_{k=m+1}^n x_{\substack{{}\\ k}}\Big\| =0\, .$
\begin{itemize}
\item[(ee)] The norm-closure of the convex hull
\begin{center}
$\text{\rm conv}\,\Big(\big\{ \sigma_{\substack{{}\\ \!\leftarrow}}^k
\big( (x_{\substack{{}\\ j}})_{\substack{{}\\ j\geq 1}}\big)
\, ;\, k\geq 0\big\}\Big)\subset l^\infty (X)$
\end{center}
contains the zero sequence.
\end{itemize}
\end{theorem}

\begin{proof} Without loss of generality we can assume that $\| x_k\|
\leq 1$ for all $k\geq 1\,$.

The proof of (e)$\,\Rightarrow\,$(ee) is immediate. Indeed, if
$\varepsilon >0$ and $n_\varepsilon\in\mathbb{N}^*$ are such that
\smallskip

\noindent\hspace{2.046 cm}$\displaystyle \Big\| \sum_{k=m+1}^n
x_{\substack{{}\\ k}}\Big\| \leq (n-m)\,\varepsilon\, ,\qquad 0\leq
m<n\, ,\; n-m\geq n_\varepsilon\, ,$
\medskip

\noindent then we have for every $n\geq n_\varepsilon\,$:
\begin{equation*}
\Big\|\,\frac 1n\,\sum_{k=1}^n\,
\sigma_{\substack{{}\\ \!\leftarrow}}^k
\big( (x_{\substack{{}\\ j}})_{\substack{{}\\ j\geq 1}}\big)\Big\|
=\Big\|\,\frac 1n\,\sum_{k=1}^n\,
\big( x_{\substack{{}\\ j+k}}\big)_{\substack{{}\\ j\geq 1}}\Big\|
=\sup_{j\geq 1}\,\Big\|\,\frac 1n\,\sum_{k=1}^n\,
x_{\substack{{}\\ k+j}}\,\Big\|\leq\varepsilon\, .
\end{equation*}

Conversely, let us assume that (ee) is satisfied and let
$\varepsilon >0$ be arbitrary. Then there exist
$p\in\mathbb{N}^*$ and $\lambda_1\, ,\,\ldots\, ,\,\lambda_p\geq 0\, ,
\,\lambda_1+ \ldots +\lambda_p=1\,$,
such that
\begin{equation}\label{hull}
\sup_{k\geq 1}\,\Big\|\sum_{j=1}^p\,\lambda_j\, x_{j+k}\,\Big\| =
\Big\|\sum_{j=1}^p\,\lambda_j\,\sigma_{\substack{{}\\ \!\leftarrow}}^j
\big( (x_{\substack{{}\\ k}})_{\substack{{}\\ k\geq 1}}\big)\Big\|
\leq\frac{\varepsilon}2\, .
\end{equation}
On the other hand, we have for every $0\leq m<n$ with $n-m\geq p\,$:
\medskip

\noindent\hspace{1.083 cm}$\displaystyle \frac 1{n-m}\sum_{k=m+1}^n
x_{\substack{{}\\ k}} -\frac 1{n-m}\sum_{k=m+1}^n\!\Big(\,
\sum_{j=1}^p \lambda_j\, x_{\substack{{}\\ j+k}}\Big) =$
\smallskip

\noindent\hspace{0.712 cm}$\displaystyle =\frac 1{n-m}\sum_{k=m+1}^n
\, \sum_{j=1}^p\, \lambda_j\,\big( x_{\substack{{}\\ k}} -
x_{\substack{{}\\ j+k}}\big) =\frac 1{n-m}\,\sum_{j=1}^p\,
\lambda_j\!\sum_{k=m+1}^n\! \big( x_{\substack{{}\\ k}} -
x_{\substack{{}\\ j+k}}\big) =$
\smallskip

\noindent\hspace{0.712 cm}$\displaystyle =\frac 1{n-m}\,\sum_{j=1}^p
\lambda_j \Big( \sum_{k=m+1}^{m+j} x_{\substack{{}\\ k}} -
\sum_{k=n+1}^{n+j} x_{\substack{{}\\ k}}\Big)\, ,$
\medskip

\noindent hence
\begin{equation}\label{difference}
\Big\|\,\frac 1{n-m}\sum_{k=m+1}^n x_{\substack{{}\\ k}} -
\frac 1{n-m}\sum_{k=m+1}^n\!\Big(\,\sum_{j=1}^p
\lambda_{\substack{{}\\ j}}\, x_{\substack{{}\\ j+k}}\Big)\Big\|
\leq\frac{2\, p}{n-m}\, .
\end{equation}
Now (\ref{hull}) and (\ref{difference}) yield
\smallskip

\noindent\hspace{1.8 cm}$\displaystyle 0\leq m<n\, ,\,
n-m\geq\frac{4\, p}\varepsilon\quad\Longrightarrow\quad
\Big\|\,\frac 1{n-m}\sum_{k=m+1}^n x_{\substack{{}\\ k}}\Big\|
\leq\varepsilon\; .$

\end{proof}

For convex shift-bounded vector sequences the statement of Theorem
\ref{m.erg.} can be strengthened:

\begin{corollary}[Mean Ergodic Theorem for convex shift-bounded
sequences]\label{pow.bdd.m.erg.}
For a convex shift-bounded sequence
$(x_{\substack{{}\\ k}})_{\substack{{}\\ k\geq 1}}$ in a
Banach space $X\,$, the following conditions are equivalent $:$
\smallskip

\noindent\hspace{0.62 cm}{\rm (e)}\hspace{3.1 cm}
$\displaystyle \lim_{\substack{0\leq m<n\\ n-m\to\infty}}\; \Big\|
\,\frac 1{n-m}\,\sum_{k=m+1}^n x_{\substack{{}\\ k}}\Big\| =0\, .$
\smallskip

\noindent\hspace{0.62 cm}{\rm (f)}\hspace{3.9 cm}
$\displaystyle \lim_{n\to\infty}\; \Big\|\,\frac 1n\,\sum_{k=1}^n
x_{\substack{{}\\ k}}\Big\| =0\, .$
\begin{itemize}
\item[(ff)] The weak closure of the convex hull
$\,\text{\rm conv}\,\big( \{ (x_{\substack{{}\\ k}}\, ;\, k\geq 1\}
\big)\subset X$ contains $0\,$.
\end{itemize}
\end{corollary}

\begin{proof}
The implications (e)$\,\Rightarrow\,$(f)$\,\Rightarrow\,$(ff)
are trivial.

Since the weak closure of $\,\text{\rm conv}\,\big( \{
(x_{\substack{{}\\ k}}\, ;\, k\geq 1\}\big)$ is equal to its norm
closure, (ff) implies that, for any $\varepsilon >0\,$, there
are $p\in\mathbb{N}^*$ and $\lambda_1\, ,\,\ldots\, ,\,\lambda_p
\geq 0\, ,\,\lambda_1+ \ldots +\lambda_p=1\,$, such that

\noindent\hspace{5 cm}$\displaystyle \Big\|\sum_{j=1}^p\,
\lambda_j\, x_j\,\Big\|\leq\varepsilon\, .$
\medskip

\noindent Using (\ref{pow.bdd}), it follows that
\medskip

\noindent\hspace{2.42 cm}$\displaystyle \Big\|\sum_{j=1}^p\,
\lambda_j\,\sigma_{\substack{{}\\ \!\leftarrow}}^j\big(
(x_{\substack{{}\\ k}})_{\substack{{}\\ k\geq 1}}\big)\Big\| =
\sup_{k\geq 1}\,\Big\|\sum_{j=1}^p\,\lambda_j\, x_{j+k}\,\Big\|
\leq c\,\varepsilon\, .$
\medskip

By the aboves (ff) implies condition (ee) in Theorem \ref{m.erg.},
hence (e).

\end{proof}

%

%
%
%

\bigskip
\subsection*{Acknowledgment}
Part of this paper was written while the author was guest at
l'Universit\'e de Lille in March 2002. He is grateful to
Professors Mustafa Mbekhta and Florian-Horia Vasilescu for
their warm hospitality.

\bigskip
\bibliographystyle{amsplain}

\end{document}